\documentclass[12pt]{article}
\usepackage{graphicx,amsfonts,amsmath,amssymb,amsthm,verbatim,float,subfig,fullpage,algorithmic,xcolor,natbib,xr,xcolor,psfrag,epsf,bbm,bm,times}
\usepackage[latin1]{inputenc}
\usepackage[plain,noend]{algorithm2e}
\usepackage[shortlabels]{enumitem}
\usepackage[hidelinks]{hyperref}

\usepackage{cleveref}
\crefname{assumption}{Assumption}{Assumptions}
\crefname{proposition}{Proposition}{Propositions}
\crefname{theorem}{Theorem}{Theorems}
\crefname{section}{Section}{Section}
\crefname{lemma}{Lemma}{Lemma}
\crefname{algorithm}{Algorithm}{Algorithms}
\crefname{appendix}{Appendix}{Appendix}
\crefname{example}{Example}{Examples}
\crefname{remark}{Remark}{Remarks}
\crefname{figure}{Figure}{Figure}
\crefname{equation}{equation}{equation}

\newtheorem{theorem}{Theorem}  
  
\newtheorem{proposition}{Proposition}

\newtheorem{assumption}{Assumption}

\allowdisplaybreaks

\def\1{\mathbf{1}}

\renewcommand{\phi}{\varphi}

\newcommand{\B}{\mathcal{B}}

\newcommand{\F}{\mathcal{F}}
\newcommand{\G}{\mathcal{G}}

\newcommand{\M}{\mathcal{M}}
\newcommand{\N}{\mathrm{Normal}} 
\renewcommand{\P}{\mathcal{P}}

\newcommand{\X}{\mathcal{X}}
\newcommand{\Y}{\mathcal{Y}}

\def\EE{E}
\def\PP{\mathsf{pr}}
\renewcommand{\O}{\mathcal{O}}
\newcommand{\RR}{\mathbb{R}}
\newcommand{\VV}{\mathbb{V}}
\newcommand{\XX}{\mathsf{X}}
\newcommand{\YY}{\mathsf{Y}}

\newcommand{\indicator}{\mathbb{I}}

\newcommand{\dd}{\text{d}}
\newcommand{\dx}{\text{d}x}

\newcommand{\initial}{\pi_0}
\newcommand{\kernel}{K}

\newcommand{\potential}{g}

\newcommand{\forward}{\mathsf{F}}
\newcommand{\forwardQ}{\mathsf{Q}}
\newcommand{\forwardQQ}{ \widetilde{\mathsf{Q}} }

\newcommand{\bra}[1]{\langle #1 \rangle}
\newcommand{\BK}[1]{ {\left( #1 \right)} }
\newcommand{\sqBK}[1]{ {\left[ #1 \right]} }
\newcommand{\curBK}[1]{ {\left\{ #1 \right\}} }

\newcommand{\norm}[1]{\left\Vert #1 \right\Vert}
\newcommand{\tnorm}[1]{{\vert\kern-0.25ex\vert\kern-0.25ex\vert #1 \vert\kern-0.25ex\vert\kern-0.25ex\vert}}

\newcommand{\vf}[1]{\textbf{\color{teal} \textsc{[verified]}}}

\def\Var{\mathrm{var}}

\newcommand{\plim}{\ensuremath{\stackrel{\mathrm{pr}}{\rightarrow}}}

\newcounter{hypA}

\newcommand{\iu}{\mathrm{i}\mkern1mu}

\def\Wbar{\overline{W}}
\def\pihat{\widehat{\pi}}

\def\b1{\mathbf{1}}
\def\eqd{\stackrel{\mathrm{D}}{=}}
\def\alphadist{\upsilon}

\def\Zhat{\widehat{Z}}
\def\gammahat{\widehat{\gamma}}
\def\bistochm{\mathcal{M}}
\def\muhat{\widehat{\mu}}

\renewcommand\mod{\, \mathrm{mod} \,}

\def\SS{\mathcal{S}}

\def\v{{\varepsilon}}

\begin{document}
\title{Particle filter efficiency under limited communication}
\author{Deborshee Sen\\ {ds2469@bath.ac.uk}  
}

\date{Department of Mathematical Sciences, University of Bath, Bath BA27AY, UK}

\maketitle

\begin{abstract}
Sequential Monte Carlo methods are typically not straightforward to implement on parallel architectures. This is because standard resampling schemes involve communication between all particles.
The $\alpha$-sequential Monte Carlo method was proposed recently as a potential solution to this which limits communication between particles. This limited communication is controlled through a sequence of stochastic matrices known as $\alpha$-matrices. We study the influence of the communication structure on the convergence and stability properties of the resulting algorithms. In particular, we quantitatively show that the mixing properties of the $\alpha$-matrices play an important role in the stability properties of the algorithm. Moreover, we prove that one can ensure good mixing properties by using randomized communication structures where each particle only communicates with a few neighboring particles. The resulting algorithms converge at the usual Monte Carlo rate. This leads to efficient versions of distributed sequential Monte Carlo.
\end{abstract}


\section{Introduction}

Hidden Markov models \citep{rabiner1986introduction}, also known as state-space models \citep{durbin2012time}, constitute a large class of numerical methods frequently used in statistics and signal processing. Examples of application areas include ecology \citep{michelot2016movehmm}, finance \citep{nystrup2017long}, medical physics \citep{ingle2015ultrasonic}, natural language processing \citep{kang2018opinion}, oceanology \citep{grecian2018understanding}, and sociology \citep{qiao2017predicting}.

A hidden Markov model with measurable state space $ (\XX,\mathcal{X})$ and observation space $ (\YY,\mathcal{Y})$ is a process $ \{ (X_t,Y_t ) \}_{t \geq 0}$, where $ \{ X_t \}_{t \geq 0} $ is a Markov chain on $\XX$, and each observation $Y_t$, valued in $\YY$, is conditionally independent of the rest of the process given $X_t$. Let $\initial$ and $\{ K_t \}_{t \geq 1}$ be respectively a probability distribution and a sequence of Markov kernels on $(\XX, \X )$, and let $\{ g_t \}_{t \geq 0}$ be a sequence of Markov kernels acting from $ (\XX,\X )$ to $(\YY, \Y )$, with $g_t (x,\cdot)$ admitting a strictly positive density -- denoted similarly by $g_t (x,y)$ -- with respect to some dominating $\sigma$-finite measure for every $t \geq 0$, which we shall assume to be the Lebesgue measure for convenience. The hidden Markov model specified by $\initial$, $\{ K_t \}_{t \geq 1}$ and $\{ g_t \}_{t \geq 0}$ is 
\begin{align} \label{eq.hmm}
\begin{aligned}
X_0 
& \sim 
\initial(\cdot), 
\\ 
X_t \mid ( X_{t-1} = x_{t-1}) 
& \sim 
K_t (x_{t-1}, \cdot) ~~
(t \geq 1), 
\\
Y_t \mid ( X_t = x_t ) 
& \sim 
g_t (x_t,\cdot) ~~ 
(t \geq 0).
\end{aligned}
\end{align}

In the sequel, we fix a sequence of observations $y = \{ y_t \}_{t \geq 0}$ and use $g_t (x)$ to denote $g_t (x, y_t)$ for $t \geq 0$.
The functions $\{g_t(\cdot)\}_{t \geq 0}$ are known as potential functions and the kernels $\{ K_t \}_{t \geq 1}$ are known as latent transition kernels. 
Let $\M(\XX)$ and $\P(\XX)$ denote the set of measures and probability measures on $(\XX,\X)$, respectively, and let $\B(\XX)$ denote the set of all real-valued measurable functions on $(\XX,\X)$ which are bounded by one in absolute value. For a measure $\pi \in \M(\XX)$ and a function $\phi \in \B(\XX)$, we define $\pi(\phi) = \int_\XX \phi(x) \pi (\dd x)$, and for a Markov kernel $K$ on $(\XX, \X)$, we define $K \phi(x) = \int_\XX \phi(x') K(x, \dd x')$. We use the notation $Y_{s:t}$ for $s \leq t$ to denote $(Y_s, \dots, Y_t)$.

We focus our attention on the predictive distribution in this article, which is the distribution of $X_T \mid Y_{0:(T-1)}$ for $T \geq 1$. The analysis developed can be straightforwardly extended to the filtering distribution, which is the distribution of $X_T \mid Y_{0:T}$. We denote the predictive distribution by $\pi_T\{X_T \mid Y_{0:(T-1)}\}$ for $T \geq 1$. Integrals of functions $\phi \in \B (\XX)$ with respect to the predictive distribution can be written as 
\begin{equation} \label{eq:normalised_target}
\pi_T (\phi) 
= 
\frac{1}{Z_T} \int_{\XX^{T+1}} \initial(\dd x_0) \prod_{t=1}^T K_t(x_{t-1}, \dd x_t) \prod_{t=0}^{T-1} g_t (x_t) \, \phi(x_T) ,
\end{equation}
where $Z_T$ is the normalisation constant, which is the marginal likelihood of the observations $Y_{0:(T-1)}$ given by $Z_T = \int_{\XX^{T+1}} \initial(\dd x_0) \prod_{t=1}^T K_t (x_{t-1}, \dd x_t) \prod_{t=0}^{T-1} g_t (x_t)$; we also define $Z_0 = 1$.
For the purpose of analysis, it is useful to also consider the unnormalised measure $\gamma_T$ defined as 
\begin{equation} \label{eq:unnormalised_target}
\gamma_T(\phi) 
=
Z_T \times \pi_T(\phi).
\end{equation}

Unfortunately, these integrals cannot be evaluated analytically except for linear Gaussian models, and Monte Carlo methods must be used instead. The bootstrap particle filter algorithm \citep{gordon1993novel} is commonly used for inference in hidden Markov models. It starts by generating $N \geq 1$ independent and identically distributed samples, termed \emph{particles}, $X_0 = \{ X_0^i \}_{i=1}^N$ from the distribution $\pi_0$. Given particles $ X_{t-1} = \{X_{t-1}^i\}_{i=1}^N$, it performs multinomial \emph{resampling} according to (unnormalised) weights $\{ g_{t-1}(X_{t-1}^i)\}_{i=1}^N$, before propagating the particles via the Markov kernel $K_t$. At each time $t \geq 0$, the bootstrap particle filter provides a particle approximation of the predictive distribution $\pi_t$ and the normalisation constant $Z_t$.

Parallel and distributed algorithms have become increasingly relevant as parallel computing architectures have become the norm rather than the exception. While there has been significant research devoted to distributed Markov chain Monte Carlo algorithms \citep{ahn2014distributed, scott2016bayes, li2017simple, heng2019unbiased, ou2021scalable}, the same has generally not been true for particle filtering. 
The resampling step of particle filters makes it difficult to parallelize. 
Two parallel implementations of the resampling step were proposed by \cite{bolic2005resampling}, and alternative schemes were investigated by \cite{miao2011algorithm, murray2012gpu, murray2016parallel}. 
\cite{verge2015parallel} provided algorithms involving resampling at two hierarchical levels, and \cite{del2017convergence} proved convergence and central limit theorem. \cite{miguez2014uniform, miguez2016proof} provided proofs of convergence for distributed particle filters relying on techniques developed in \cite{bolic2005resampling}; however, these assumed that a certain notion of weight degeneracy does not occur. We prove in this article that weight degeneracy can be avoided by suitably choosing network architectures of distributed particle filters. \cite{heine2020parallelizing} designed stable-in-time distributed sequential Monte Carlo algorithms with limited interactions; however, these converge at a slower rate than the standard Monte Carlo rate.

The $\alpha$-sequential Monte Carlo algorithm \citep{whiteley2016role} was proposed recently as a general method for distributed sequential Monte Carlo. This is a generalisation of the bootstrap particle filter that can be implemented on parallel architectures. 
This is achieved by allowing particles to interact with only a small subset of other particles in the resampling step, and is formalized through a sequence of stochastic matrices. These are referred to as {\it connectivity matrices} in the sequel since they describe how particles are connected to each other. 
It has been shown that certain ``local exchange'' communication structures do not lead to stable algorithms \citep{heine2017fluctuations}, and sophisticated adaptive mechanisms have been designed for ensuring stability \citep{lee2016forest, heine2020parallelizing}. However, a general understanding of the influence of the communication structure on the stability properties of the algorithm is lacking. 

In this article, we relate the stability properties of the $\alpha$-sequential Monte Carlo algorithm to the connectivity and mixing properties of the communication structures described by the connectivity matrices. In particular, we show that it is possible to design $\alpha$-sequential Monte Carlo algorithms with time-uniform convergence at the standard Monte Carlo rate of $N^{-1/2}$ without the degree of the interaction graph growing with the number of particles $N$. Computer code for numerical experiments in this article can be found online at \url{https://github.com/deborsheesen/alphaSMC}.

\section{\texorpdfstring{$\alpha$}--Sequential Monte Carlo} \label{sec.alphaSMC}

\subsection{Algorithm description}

The $\alpha$-sequential Monte Carlo algorithm with $N \geq 1$ particles relies on a sequence of (possibly random) matrices $\{\alpha_t\}_{t \geq 0}$, where each $\alpha_t = ( \alpha_t^{ij} )_{i, j = 1}^N \in \RR^{N,N}$ is a stochastic matrix; for any time index $t \geq 0$ and particle index $i = 1, \dots, N$, we have $\sum_{j=1}^N \alpha_t^{ij} = 1$. The $\alpha$-sequential Monte Carlo algorithm simulates a sequence $\{ X_t; \, t \geq 0 \}$, where, for each time index $t \geq 0$, we have $X_t = \{ X_{t}^{i} \, ; \, i = 1, \dots, N \}$, and $X^i_t \in \XX$ is the location of the $i$-th particle at time index $t \geq 0$. The particle approximation $\pihat_t^N$ of $\pi_t$ produced by the $\alpha$-sequential Monte Carlo algorithm is given by $\pihat_t^N = \sum_{i=1}^N \Wbar_t^i \, \delta_{X_{t}^{i}}$, 
where $\Wbar_t = (\Wbar^1_t, \dots, \Wbar^N_t) \in \mathcal{P}(N)$ denotes the vector of normalised weights with $\P(N) = \{x \in \RR^N_+ \; : \; \sum_{i=1}^N x_i = 1\}$
being the $N$-dimensional probability simplex. 
We have also defined $\Wbar^i_t = W^i_t / (\sum_{j=1}^N W^j_t)$ as the normalised weights. The unnormalised weights $W_t = (W^1_t, \dots, W^N_t) \in \RR_+^N$ are recursively defined as follows. At time index $t=0$, the weights are all initialized to one, that is, $W^i_t = 1$ $(i=1, \dots, N)$. For $t \geq 1$, the weights are recursively defined as
\begin{equation} \label{eq:alphaSMC_weights}
W_{t}^i = \sum_{j=1}^N \alpha_{t-1}^{ij} \, W_{t-1}^{j} \, g_{t-1}(X_{t-1}^{j}) ~~ (i=1, \dots, N).
\end{equation}

The $\alpha$-sequential Monte Carlo algorithm also produces a particle approximation of the unnormalised measure $\gamma_t$ and the normalisation constant $Z_t$ as $\gammahat_t^N = (1/N) \sum_{i=1}^N W_t^i \, \delta_{X_{t}^{i}}$ and $\Zhat^N_t = \gammahat_t^N(1) 
= (1/N) \sum_{i=1}^N W_t^i$.
%
%
The particle equivalent of \cref{eq:unnormalised_target} is $\gammahat_t^N = \Zhat_t^N \times \pihat_t^N$, which states that the estimate of the unnormalised measure can be decomposed into the product of estimates for the normalised measure and the normalisation constant; this is the same as for the bootstrap particle filter.

The particles are initialised as follows. At time index $t=0$, particles $X_0^i \in \XX$ are simulated as being independent and identically distributed from the initial distribution $\pi_0$. We define $\F_{t-1}$ to be the $\sigma$-algebra generated by all the particles up to and including time $(t-1)$, that is, $X_{0:(t-1)}$, and all the connectivity matrices up to and including time $(t-1)$, that is, $\alpha_{0:(t-1)}$. We also define the notations $\EE_{t}(\cdot) = \EE(\cdot \mid \F_t)$ and $\Var_t(\cdot) = \Var (\cdot \mid \F_t)$ for convenience,
which are the conditional mean and variance conditioned upon the state of the system up to and including time $t$; these will typically be used in the context of events happening after time $t$.
At time index $t \geq 1$ and conditionally upon $\F_{t-1}$, the particles $\{ X_t^i \}_{i=1}^N$ are simulated independently, with
\begin{equation*}
\mathbb{P}\BK{X_t^{i} \in \dx \mid \F_{t-1}} = 
\frac{1}{W^i_t} \, \sum_{j=1}^N \alpha_{t-1}^{ij} \, W_{t-1}^{j} \, g_{t-1}(X_{t-1}^j) \, K_t(X_{t-1}^j, \dx).
\end{equation*}
The $\alpha$-sequential Monte Carlo algorithm is summarised in \cref{algo:alpha_SMC}. Throughout this text, we assume that the connectivity matrices $\{\alpha_t\}_{t \geq 0}$ can all be generated at the start of the algorithm. In other words, we do not consider adaptive schemes for constructing the connectivity matrices, as for example is explored by \cite{liu1995blind, whiteley2016role, lee2016forest}.

%
%
%
%
%

\begin{algorithm}
\caption{$\alpha$-sequential Monte Carlo algorithm \citep{whiteley2016role}.} 
\label{algo:alpha_SMC}

\textbf{Input:} Connectivity matrices $\{\alpha_t\}_{t \geq 0}$, potential functions $\{g_t\}_{t \geq 0}$, Markov kernels $\{K_t\}_{t \geq 1}$ and initial distribution $\pi_0$.

\begin{algorithmic}[1] 

\FOR{$(i = 1, \dots, N)$}

\STATE 
Set $W_0^i = 1$.

\STATE 
Sample $X_0^i \sim \initial$ independently.

\ENDFOR 

\FOR{$t \geq 1$}

\FOR{$(i = 1, \dots, N)$}

\STATE 
Set $W_{t}^i = \sum_{j=1}^N \alpha_{t-1}^{ij} \, W_{t-1}^{j} \, g_{t-1} ( X_{t-1}^{j} )$.

\STATE 
Sample 
\begin{equation*}
X_{t}^{i} \mid \F_{t-1} 
\sim 
\frac{1}{W_{t}^i} \sum_{j=1}^N \alpha_{t-1}^{ij} \, W_{t-1}^{j} \, g_{t-1} (X_{t-1}^{j} ) \, K_t ( X_{t-1}^{j}, \cdot ) ~ \text{independently.}
\end{equation*} 

\ENDFOR 

\ENDFOR 

\end{algorithmic}

\textbf{Output:} Weighted particle system $\{ (X_{t}^{i}, W_{t}^i) \, ; \, i = 1, \dots, N, ~ t \geq 1 \}$.

\end{algorithm}

\cite{zhang2020performance} have implemented a distributed resampling technique using a message passing interface for a scheme that is similar to the local exchange scheme analysed by \cite{heine2017fluctuations}, and have reported computational gains from doing so.

\subsection{Basic Properties}

The predictive probability distributions $\{ \pi_t \}_{t \geq 0}$ defined by the state-space model \eqref{eq:normalised_target} satisfy $\pi_t = \forward_t \pi_{t-1}$, where the mapping $\forward_t : \P(\XX) \rightarrow \P(\XX)$ associates to any probability measure $\pi \in \P(\XX)$ the probability measure $\forward_t \pi$ that acts on functions $\varphi\in \B(\XX)$ as
\begin{equation*} 
\forward_t \pi \, (\varphi) 
=
\frac{ \pi ( g_{t-1} K_t \varphi ) }{ \pi (g_{t-1}) }, \quad t \geq 1.
\end{equation*}
For two time indices $0 \leq s \leq t$, set $\forward_{s,t} = \forward_t \circ \cdots \circ \forward_{s+1}$, with the convention that $\forward_{t,t}$ is the identity mapping, so that we have $\pi_t = \forward_{s,t} \, \pi_s$. Similarly, the unnormalised measures $\{ \gamma_t \}_{t \geq 0}$ satisfy $\gamma_t(\phi) = \gamma_s(\forwardQ_{s,t} \, \phi)$, where $\forwardQ_{s,t} = \forwardQ_{s+1} \circ \cdots \circ \forwardQ_t$ and the operator $\forwardQ_t$ acts on a test function $\phi \in \B(\XX)$ as $\forwardQ_t \, \phi = g_{t-1} \, K_t \phi$, $t \geq 1$.

As noted in \cite{whiteley2016role}, if the connectivity matrices $\{ \alpha_{t} \}_{t \geq 0}$ keep (almost surely) the uniform distribution on $\{1, \dots, N\}$ invariant, that is, $\b1 \alpha_t = \b1$, where $\b1 = (1, \dots, 1) \in \RR^N$ is the $N$-dimensional vector of ones. The definition \eqref{eq:alphaSMC_weights} of the weights shows that the particle approximations $\gammahat_t^N$ are such that for any test function $\phi \in \B(\XX)$,
\begin{equation} \label{eq.gamma.hat.flow}
\EE_{t-1}\curBK{\gammahat^N_t(\phi)}
=
\EE_{t-1}\curBK{W^1_{t} \phi(X^1_t)}
=
\frac{1}{N} \sum_{j=1}^N W^j_{t-1} \forwardQ_t \phi(X^j_{t-1}) = \gammahat^N_{t-1}(\forwardQ_t \phi).
\end{equation}
Consequently, iterating \cref{eq.gamma.hat.flow} shows that the particle approximation $\gammahat^N_t(\phi)$ is unbiased: $\EE\{\gammahat^N_t(\phi)\} = \EE\{\gammahat^N_0( \forwardQ_{0,t} \phi)\} = \gamma_0(\forwardQ_{0,t} \phi) = \gamma_t(\phi)$.
Since $\Zhat^N_t = \gamma^N_t(1)$, it also follows that $\EE(\Zhat^N_t) = Z_t$. This lack-of-bias property allows the $\alpha$-sequential Monte Carlo approach to be straightforwardly leveraged within other Monte Carlo schemes such as the pseudo-marginal Monte Carlo approach \citep{andrieu2009pseudo}, particle Markov chain Monte Carlo methods \citep{andrieu2010particle}, and advanced sequential Monte Carlo methods \citep{chopin2013smc2}.

%
%
%
%

\section{Time-uniform stability of \texorpdfstring{$\alpha$}{}-sequential Monte Carlo} 
\label{sec.non-asymptotic_analysis}

\subsection{Mixing of connectivity matrices}

In this section, we assume that there exists a fixed bi-stochastic matrix $\alpha \in \RR_+^{N,N}$ such that $\alpha_t = \alpha$ for all $t \geq 0$. Under the assumption that the uniform distribution on $\{1, \dots, N \}$ is the unique invariant distribution of $\alpha$, we relate the stability properties of the $\alpha$-sequential Monte Carlo algorithm to the mixing properties of the connectivity matrix $\alpha$. We define the {\it mixing constant} $\lambda(\alpha)$ of the connectivity matrix $\alpha$ as
\begin{equation} \label{eq.constant.lambda}
\lambda(\alpha) 
= 
\sup_{v \in \mathsf{B}^0_1} \; \| \alpha v \| < 1,
\end{equation}
where $\| \cdot \|$ denotes the Euclidean norm and $\mathsf{B}^0_1 = \{v \in \RR^N \; : \; \|v\| = 1 \; \textrm{and} \; \bra{v,\b1} = 0\}$ is the compact set of unit vectors that are orthogonal to the vector $\b1$. The quantity $\lambda(\alpha) \geq 0$ is the smallest constant such that for any vector $v \in \RR^N$, we have
\begin{equation} \label{eq.mixing.alpha}
\norm{\alpha^k v - \BK{\sum_{i=1}^N \frac{v_i}{N}} \times \b1}
\leq 
\lambda(\alpha)^k \norm{v - \BK{\sum_{i=1}^N \frac{v_i}{N}} \times \b1}.
\end{equation}

If the Markov transition matrix $\alpha$ is reversible with respect to the uniform distribution on $\{1, \dots, N \}$, that is, $\alpha$ is symmetric, the quantity $\lambda(\alpha)$ equals the absolute value of the second largest (in absolute value) eigenvalue of $\alpha$: $\lambda(\alpha) = \max_{k \in \{2, \dots, N\}} \; |\lambda_k|$, 
where $1=\lambda_1 \geq \cdots \geq \lambda_N > -1$ is the spectrum of $\alpha$. In other words, in the reversible case, $\lambda(\alpha)$ can also be expressed as one minus the absolute spectral gap of the matrix $\alpha$. In the case where $w \in \RR^N$ is a probability vector, that is, $w \in \P(N)$, \cref{eq.mixing.alpha} can be reformulated as
\begin{equation} \label{eq.fundamental.bound}
\|\alpha w\|^2 
\leq \frac{1-\lambda^2(\alpha)}{N} + \lambda^2(\alpha) \|w \|^2.
\end{equation}
This is the key inequality that we will use to establish the stability properties of the $\alpha$-sequential Monte Carlo algorithm.

\subsection{Stability} 
To measure the discrepancy between two (possibly random) probability measures $\mu$ and $\nu$, consider the norm
\begin{equation*} 
\tnorm{\mu - \nu}^2 = \sup \curBK{
\EE \left [ \curBK{ \mu(\varphi) - \nu(\varphi) }^2 \right ] \; : \; \varphi \in \B (\XX) }.
\end{equation*}
We assume in this section that the potential functions $\{g_t\}_{t \geq 0}$ and latent transition kernels $\{K_t\}_{t \geq 1}$ of the state-space model \eqref{eq.hmm} are uniformly bounded in time; this is standard when studying the stability properties of particle filters \citep{del2001stability, whiteley2016role}. In other words, we make the following \cref{assumption_1}.
\begin{assumption} \label{assumption_1}
There exist constants $\kappa_K > 1$ and $\kappa_g > 1$ such that $\kappa_K^{-1} \leq K_t \leq \kappa_K$ $(t \geq 1)$ and $\kappa_g^{-1} \leq g_t \leq \kappa_g$ $(t \geq 0)$.
\end{assumption}
The main result of this section is that under \cref{assumption_1} and as soon as the absolute spectral gap of the matrix $\alpha$ is large enough, the discrepancy between the particle approximation $\pihat_t^N$ and its limiting value $\pi_t$ can be uniformly bounded in time: 
\begin{equation*}
\sup_{t \geq 0} \tnorm{\pihat^N_t - \pi_t} 
\leq
\textrm{Cst} \times N^{-1/2}.
\end{equation*}
In other words, the particle approximation $\pihat^N_t$ converges to the true predictive distribution $\pi_t$ at the usual Monte Carlo rate, and this convergence can be controlled uniformly in time. This is formalised in \cref{thm:non-asymptotic_stability}, which is proved in \cref{proof:non-asymptotic_stability}.

%
%
\begin{theorem}[Uniform stability] \label{thm:non-asymptotic_stability}
Suppose that the state-space model \eqref{eq.hmm} satisfies Assumption \ref{assumption_1}. Consider the $\alpha$-sequential Monte Carlo algorithm with $N$ particles and a constant bi-stochastic connectivity matrix $\alpha \in \RR_+^{N,N}$ such that
\begin{itemize}
\item 
the uniform distribution on $\{1, \dots, N \}$ is the unique invariant distribution of $\alpha$, and

\item 
the mixing constant $\lambda(\alpha)$ defined in \cref{eq.constant.lambda} satisfies $\lambda(\alpha) < \kappa_g^{-2}$.
\end{itemize}
Then the following uniform bound for the $N$-particle approximations $\pihat_t$ holds:
\begin{equation}\label{eq.stability.bound}
N \times \tnorm{\pihat^N_t - \pi_t}^2
\; \leq \;
\frac{D}{1-\rho} \times \frac{\kappa_g^4 \{1 - \lambda^2(\alpha)\}}{1 - \kappa_g^4 \lambda^2(\alpha)} 
\end{equation}
for constants $D>0$, $\kappa_g>1$ and $0<\rho<1$ that depend only on the state-space-model \eqref{eq.hmm}.
\end{theorem}
The bootstrap particle filter corresponds to the case where $\lambda(\alpha) = 0$, and in that case one obtains that $N \times \tnorm{\pihat^N_t - \pi_t}^2 \; \leq \; D \kappa_g^4/(1-\rho) = C_{\textrm{bootstrap}}$.
In the case of fast mixing connectivity matrices, that is, $\lambda(\alpha) \ll 1$, expanding the right-hand-side of \cref{eq.stability.bound} in powers of $\lambda(\alpha)$ yields that
\begin{equation} \label{eq.communication.costs.spectral}
N \times \tnorm{\pihat^N_t - \pi_t}^2
\; \leq \;
C_{\textrm{bootstrap}} 
+ 
\textrm{Const} \times \lambda^2(\alpha) + \mathcal{O}\{\lambda^4(\alpha)\} .
\end{equation}
In other words, when compared to the bootstrap particle filter, the use of a connectivity matrix $\alpha$ with limited communication incurs a cost of leading order $\lambda^2(\alpha)$. In \cref{sec.cost.sparse}, we discuss another situation leading to similar conclusions.

\section{Randomized connectivity matrices} \label{sec:randomized_connections}

\subsection{Setting and basic properties}
\label{sec.random.basic.prop}

We extend the analysis of the previous section to randomized connectivity structures and obtain a central limit theorem.
To this end, let $\bistochm_N$ be the set of all $N \times N$ symmetric stochastic matrices, and consider a distribution $\alphadist_N$ on $\bistochm_N$ such that 
\begin{equation} \label{eq.permute.property}
\alpha \sim \alphadist_N 
\implies 
P \alpha P^{-1} \sim \alphadist_N 
~~ \text{for any} ~ N \times N ~\text{permutation matrix} ~ P.
\end{equation}
The operation $P \alpha P^{-1}$ corresponds to permuting the nodes of the graph associated with the matrix $\alpha$. For $\alpha \in \bistochm_N$, let $\SS(\alpha)$ be the set of all permutations of the nodes of the graph associated with $\alpha$. The distribution $\alphadist_N$ is uniform over the set $\SS(\alpha)$ for every $\alpha$.
Moreover, the mixing constant \eqref{eq.constant.lambda} is the same for every matrix in the set $\SS(\alpha)$, and this is therefore a generalisation of the setting considered in \cref{sec.non-asymptotic_analysis}.
Common examples of this framework include the bootstrap particle filter (which corresponds to $\alphadist_N$ placing mass one on the matrix $\b1 \b1^T/N$) and importance sampling (which corresponds to $\alphadist_N$ placing mass one on the identity matrix). We shall consider another such setting in \cref{sec.random.permutations} with limited connections. We study the asymptotic behaviour of the $\alpha$-sequential Monte Carlo algorithm under this setting. 

In order to keep the analysis simple, we assume in this section that the potential functions $\{g_t\}_{t \geq 0}$ are uniformly bounded: there exists a constant $\kappa_g > 0$ such that, for any time index $t \geq 0$ and $x \in \XX$, we have $0< g_t(x) \leq \kappa_g$; we note that this is weaker than \cref{assumption_1}. 
We also make the following assumption.
\begin{assumption} \label{ass.cross.terms}
The distribution $\alphadist_N$ is such that for $\alpha \sim \alphadist_N$, $\EE(\alpha^{ij} \alpha^{ik}) = \O(N^{-2})$ for all $i \neq j \neq k$. In particular, there exists $0 \leq c_3 < \infty$ such that $\EE(\alpha^{ij} \alpha^{ik}) \leq c_3 N^{-2}$ for $N$ large. 
\end{assumption}
\cref{ass.cross.terms} is clearly satisfied for the bootstrap particle filter and for sequential importance sampling.

We prove consistency and a central limit theorem for the normalised measures $\pihat_t^N$ and unnormalised measures $\gammahat_t^N$ under this setting. 
In the proof of the central limit theorem, we will need to consider a further sequence of unnormalised measures defined as $\muhat^N_t = 
(1/N) \sum_{i=1}^N (W^i_t)^2 \delta_{X^i_t}$.
Define an operator $\forwardQQ_t$ that acts on a test function $\phi \in \B(\XX)$ as $\forwardQQ_t \phi \; = \; g_{t-1}^2 K_t \phi$. We show in \cref{sec.randomised.results} that under \cref{ass.cross.terms}, as $N \to \infty$, the unnormalised measures $\muhat^N_t$ converge to the measure $\mu_t$ defined as $\mu_0 = \pi_0$, and, for a test function $\phi \in \B(\XX)$,
\begin{align} \label{eq.flow.mu.basic}
\begin{aligned}
\mu_t(\phi) 
& =
\mu_{t-1}(\forwardQQ_t \phi) 
\times
\lim_{N \to \infty} \sqBK{\EE \{(\alpha^{11})^2\} + N \EE\{(\alpha^{12})^2\}}
\\
& \quad +
Z_t^2 \pi_t(\phi) 
\times 
\lim_{N \to \infty} \curBK{ N^2 \EE(\alpha^{12} \alpha^{13}) + 2N \EE(\alpha^{11}\alpha^{12} )};
\end{aligned} 
\end{align}
we have implicitly assumed that the limits on the right hand side of \cref{eq.flow.mu.basic} exist. This is true for the bootstrap particle filter and sequential importance sampling, and more generally is true for the settings we consider in \cref{sec.statistical.tradeoff}.
Moreover, \cref{ass.cross.terms} and \cref{prop.basic} of \cref{sec.random.connect.setup} ensure that the right hand side of the previous equation is finite. 
We shall exploit \cref{eq.flow.mu.basic} to study the asymptotic behaviour of $\alpha$-sequential Monte Carlo with sparse connections in \cref{sec.cost.sparse}.

\subsection{Consistency and central limit theorem} \label{sec.randomised.results}

%
%

We first establish that the particle approximations $\pihat^N_t$, $\gammahat^N_t$, and $\muhat^N_t$ are consistent.

\begin{theorem}[Consistency] \label{thm.consistency}
Assume that the potential functions satisfy $0< g_t(x) \leq \kappa_g$, and suppose also that \cref{ass.cross.terms} holds. 
For any test function $\phi \in \B(\XX)$, as $N \to \infty$, the particle approximations $\pihat^N_t(\phi)$, $\gammahat^N_t(\phi)$ and $\muhat^N_t(\phi)$ converge in probability to $\pi_t(\phi)$, $\gamma_t(\phi)$, and $\mu_t(\phi)$, respectively.
\end{theorem}

\cref{thm.consistency} is proved in \cref{proof:consistency}.
Consistency of the particle approximations $\pihat^N_t(\phi)$ and $\gammahat^N_t(\phi)$ was established in \cite{whiteley2016role} under an asymptotic negligibility condition, which is automatically satisfied when the $\alpha$ matrices are bi-stochastic; we nonetheless include a straightforward proof for the sake of being a self-contained article. The consistency of $\muhat^N_t(\phi)$ is a more involved proof and is novel in our work.

%
%

We next show a central limit theorem for the particle approximations $\pihat^N_t$ and $\gammahat^N_t$, which is proved in \cref{proof:CLT}. 

%
%
\begin{theorem}[Central limit theorem] \label{thm.central limit theorem}
Assume that the potential functions satisfy $0< g_t(x) \leq \kappa_g$, and suppose also that \cref{ass.cross.terms} holds. For any bounded test function $\phi \in \B(\XX)$, the re-normalised quantities $N^{1/2} \{\gammahat^N_t(\phi) - \gamma_t(\phi)\}$ and $N^{1/2} \{\pihat^N_t(\phi) - \pi_t(\phi)\}$ converge in laws to centred Gaussian distributions with variances $\VV^\gamma_t(\phi)$ and $\VV^\pi_t(\phi)$, respectively, where the variances satisfy the following recursions:
\begin{align} \label{eq.asymp.var}
\begin{aligned}
\VV^\gamma_t(\phi) 
& = 
\VV^\gamma_{t-1}(\forwardQ_t \phi) + \mu_t(\phi^2) - Z^2_t \pi_{t}(\phi)^2,
\\ 
\VV^\pi_t(\phi)
& = 
\frac{ \VV^\pi_{t-1}(\forwardQ_t \overline{\phi})}{\pi_{t-1}(g_{t-1})^2} + \mu_t (\overline{\phi}^2),
\end{aligned}
\end{align}
where $\overline{\phi}_t = \phi - \pi_t(\phi)$.
\end{theorem}

\cref{thm.central limit theorem} provides a way to quantify the trade-off (relative to the bootstrap particle filter) in using $\alpha$-sequential Monte Carlo under different settings as measured by its asymptotic variance. It is worth stressing that the terms $\mu_t(\varphi^2)$ and $\mu_t(\overline{\varphi}^2)$ depend on the choice of the $\alpha$-matrices used. We discuss this in more detail in \cref{sec.statistical.tradeoff}. In particular, we consider a setting in which particles are connected to a few other particles at each time and study the effect of the number of connections on the asymptotic variances.

%
%

\section{Statistical tradeoffs} \label{sec.statistical.tradeoff}

\subsection{Permutations of a random walk on \texorpdfstring{$d$}{}-regular graph}
\label{sec.random.permutations}

We describe and analyse a version of $\alpha$-sequential Monte Carlo with sparse connections that falls into the setting considered in \cref{sec.random.basic.prop}.
Consider an undirected $d$-regular graph $\G_N$ with $N$ vertices. Let $A$ be the stochastic matrix corresponding to a random walk on $\G_N$. In other words, $A^{ij} = d^{-1}$ if nodes $i$ and $j$ have a vertex connecting them, and zero otherwise. Let $\P_{{\rm permute}}$ be the uniform distribution over the set of all permutations of $\{1, \dots, N\}$, and let $\alphadist_N$ be a distribution over $\bistochm_N$ specified by $P A P^{-1}$ for $P \sim \P_{{\rm permute}}$. The operation $P A P^{-1}$ re-indexes $A$ by the permutation $\sigma$ of the indices. We consider the case where the graph $\G_N$ corresponds to a random $d$-regular graph without self connections (that is, no node is connected to itself). In this case, $\lim_{N \to \infty} \EE(\alpha^{ii}) = 0$, $\lim_{N \to \infty} N \EE\{(\alpha^{ij})^2\} = 1/d$, $\lim_{N \to \infty} 2N \EE(\alpha^{ii} \alpha^{ij}) = 0$, and $\lim_{N \to \infty} N^2 \EE(\alpha^{ij} \alpha^{ik}) = (d-1)/d$. By \cref{eq.flow.mu.basic}, this implies 
\begin{equation} \label{eq.flow.mu} 
\mu_t(\phi) 
= 
\frac1d \mu_{t-1}(\forwardQQ_t \phi) + \frac{d-1}{d} Z_t^2 \pi_t(\phi).
\end{equation}
This is used to analyze the asymptotic variances \eqref{eq.asymp.var} of $\alpha$-sequential Monte Carlo in the next section and compare them to those of the bootstrap particle filter.

\subsection{Cost under sparse connections} 
\label{sec.cost.sparse}

We leverage the central limit theorem to analyze the influence of the number of connections $d$ on the performance of the $\alpha$-sequential Monte Carlo algorithm. 
Iterating \cref{eq.flow.mu} immediately shows that $\mu_t(\phi) = Z^2_t \pi_t(\phi) + \sum_{k=1}^t \beta_{t,k}(\phi)/d^k$
for some coefficients $\{\beta_{t,k}(\phi)\}_{k=1}^t$ that depend on the test function $\phi$ and the state-space model \eqref{eq.hmm}, but not on the connectivity $d$. It then follows from Theorem \ref{thm.central limit theorem} that the asymptotic variance can be expanded as
\begin{equation*}
\VV^\gamma_t(\phi) 
= 
\curBK{ Z^2_0 \Var_{\pi_0}(\forwardQ_{0,t} \phi) + Z^2_1 \Var_{\pi_1}(\forwardQ_{1,t} \phi) + \cdots + Z^2_t \Var_{\pi_t}(\phi)}
+ \sum_{k=1}^t \frac{\widetilde{\beta}_{t,k}(\phi)}{d^k}
\end{equation*}
for some coefficients $\{\widetilde{\beta}_{t,k}(\phi)\}_{k=1}^t$ that depend on the test function $\phi$ and the state-space model \eqref{eq.hmm}, but not on the connectivity $d$; here $\Var_{\pi_s}$ denotes the variance under $\pi_s$. Not surprisingly, since the limit $d \to \infty$ corresponds to the bootstrap particle filter, the first term on the right-hand side of the previous equation equals exactly the asymptotic variance obtained from a standard bootstrap particle filter \citep{chopin2004central}. In other words,
\begin{equation} \label{eq.cost.alpha.sequential Monte Carlo}
\VV^\gamma_t(\phi) 
= 
\VV^{\textrm{bootstrap}}_t(\phi)
+ \sum_{k=1}^t \frac{ \widetilde{\beta}_{t,k}(\phi) }{d^k}
\approx
\VV^{\textrm{bootstrap}}_t(\phi)
+ \frac{ \widetilde{\beta}_{t,1}(\phi) }{d},
\end{equation}
where $\VV^{\textrm{bootstrap}}_t(\phi)$ denotes the asymptotic variance of the bootstrap particle filter.

\begin{figure}[!ht]
\centering
\includegraphics[width=0.85\textwidth]{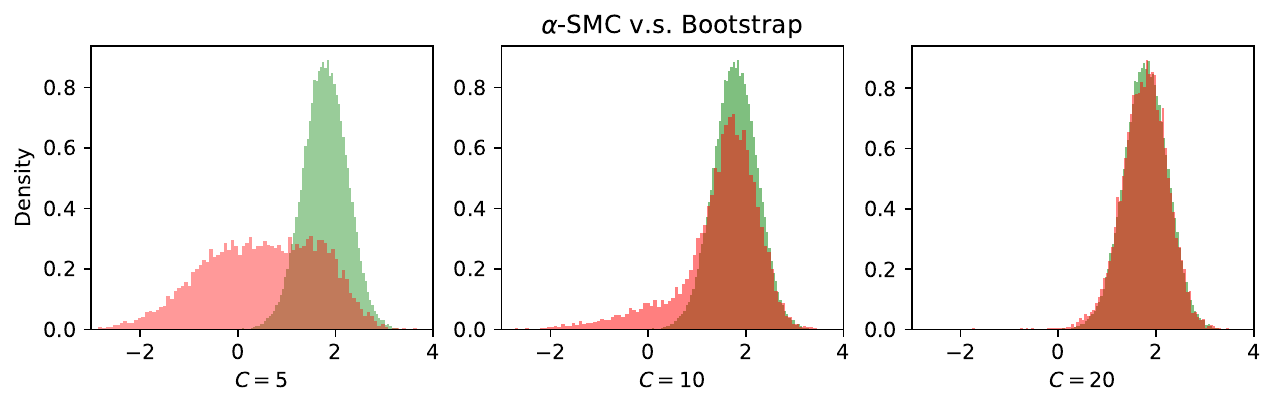}
\caption{Dynamics: $X_{t+1} = \beta X_{t} + \sqrt{1-\beta^2} \xi_t$ with $\beta=0.9$ and $\xi_t \sim \N(0,1)$. Potentials: $g_t(x) = 0.1 + 10 \times \indicator (|x-2|<0.1)$ for $t \geq 0$. We plot the accuracy of the estimation of $\pi_{T=6}(\phi)$ with $\phi(x) = x^2$.
{\bf Red:} Particle densities generated by the $\alpha$-sequential Monte Carlo with $N=10^4$ particles. 
{\bf Green:} Particle densities generated by a bootstrap particle filter.}
\label{fig.unweighted}
\end{figure}

It is interesting to note that, in general, the first coefficient $\widetilde{\beta}_{t,1}(\phi)$ can be either positive or negative. In other words, there are situations where the estimates obtained from $\alpha$-sequential Monte Carlo are statistically more efficient that those obtained from the bootstrap particle filter: $\widetilde{\beta}_{t,1}(\phi)<0$. At a heuristic level, this may be explained as follows. When using $\alpha$-sequential Monte Carlo, the propagation of information between particles is typically worse than that for the bootstrap particle filter. For example, if the distribution $\pi_t$ is more concentrated than the initial distribution $\pi_0$, it is typically the case that the distributional estimates obtained from an $\alpha$-sequential Monte Carlo with low value of $d$ will have thicker tails than the one obtained from the bootstrap particle filter (\cref{fig.unweighted}). In these situations, the $\alpha$-sequential Monte Carlo estimates of tail events of $\pi_t$
can have lower variance than the one obtained from the bootstrap particle filter.

As a concrete example, one can show that when $\pi_0$ is a standard real Gaussian distribution, $g_0(x) = 0.1 + 100 \times \indicator (|x| < 0.1)$, $\phi(x) = \indicator (|x|>1)$, and $K_1(x, \dd y) = \delta_x(\dd y)$, the $\alpha$-sequential Monte Carlo estimates of $\gamma_1(\phi) = \pi_0(g_0 \phi)$ with $d=2$ have an asymptotic variance that is roughly half as large as the one obtained from the bootstrap particle filter.

In most more realistic scenarios where particle filters are routinely used (for example, tracking of partial and/or noisy dynamical systems), though, we have indeed observed that $\alpha$-sequential Monte Carlo estimates have a higher variance than the estimates obtained from the bootstrap particle filter. Equation \eqref{eq.cost.alpha.sequential Monte Carlo} shows that there is typically a cost of order $\O (d^{-1})$ additional variance for using $\alpha$-sequential Monte Carlo instead of the bootstrap particle filter. This is demonstrated numerically in \cref{expts.alphaSMC_v_bootstrap} of \cref{sec.alphavsbootstrap}.

Connecting back to the setting considered in \cref{sec.non-asymptotic_analysis} (which considers a fixed $\alpha$ matrix), this result is in the same spirit as the bound \eqref{eq.communication.costs.spectral} that showed that there was a cost of order $\lambda(\alpha)^2$ (when controlling $N \times \tnorm{\pihat_t^N - \pi_t}^2$) when $\alpha$-sequential Monte Carlo is used instead of the bootstrap particle filter. To see the connection, consider the connectivity matrix $\alpha \in \RR_+^{N,N}$ to be equal to the Markov transition matrix of the random walk on an undirected graph $\G_N$ that is chosen uniformly at random among all the $d$-regular graphs on $N$ vertices. Any such connectivity matrix $\alpha$ is bi-stochastic, so $\lambda(\alpha)$ equals one minus the absolute spectral gap of $\alpha$: the Alon-Friedman theorem \citep{alon1986eigenvalues,friedman2008proof} states that
\begin{equation} \label{eq.alon.friedman}
\lambda(\alpha) \; \plim \; \frac{2 \sqrt{d-1}}{d} \quad \text{as} ~ N \to \infty.
\end{equation}

In other words, for such graphs and for a fixed connectivity $d \geq 2$, the mixing constant $\lambda(\alpha)$ does not deteriorate as $N \to \infty$; this is demonstrated numerically in \cref{sec.randomalpha}. If the connectivity matrix $\alpha$ was chosen this way,
for large $N$ we would observe that $\lambda^2(\alpha) = \O(d^{-1})$. \cref{thm:non-asymptotic_stability} thus shows that, under regularity assumptions on the state-space model, in order to obtain an $\alpha$-sequential Monte Carlo algorithm that is stable, one does not need to increase the number of connections $d \geq 3$ with the total number of particles $N$ as long as $d$ is large enough. 
We conjecture that a similar result holds for randomised connectivity matrices as well.
Note that if $\G_N$ is the undirected graph on $\{1, \dots, N \}$ where the vertex $i$ is connected to each vertex $j \in \{i\pm 1, \dots, i \pm \lfloor d/2 \rfloor\} \mod N$, the mixing constant $\lambda(\alpha)$ converges to one as $N \to \infty$, ultimately leading to poor performances. This is a variation of the local exchange mechanism considered in \cite{heine2017fluctuations}, where the authors indeed show that one cannot expect such an algorithm to converge uniformly at rate $N^{-1/2}$.

%
%

\section{Numerical examples} \label{sec.numerics}

\subsection{Spectral gap of random \texorpdfstring{$\alpha$}{}-matrices} \label{sec.randomalpha}

\begin{figure}[!ht]
\centering
\includegraphics[width=0.75\textwidth]{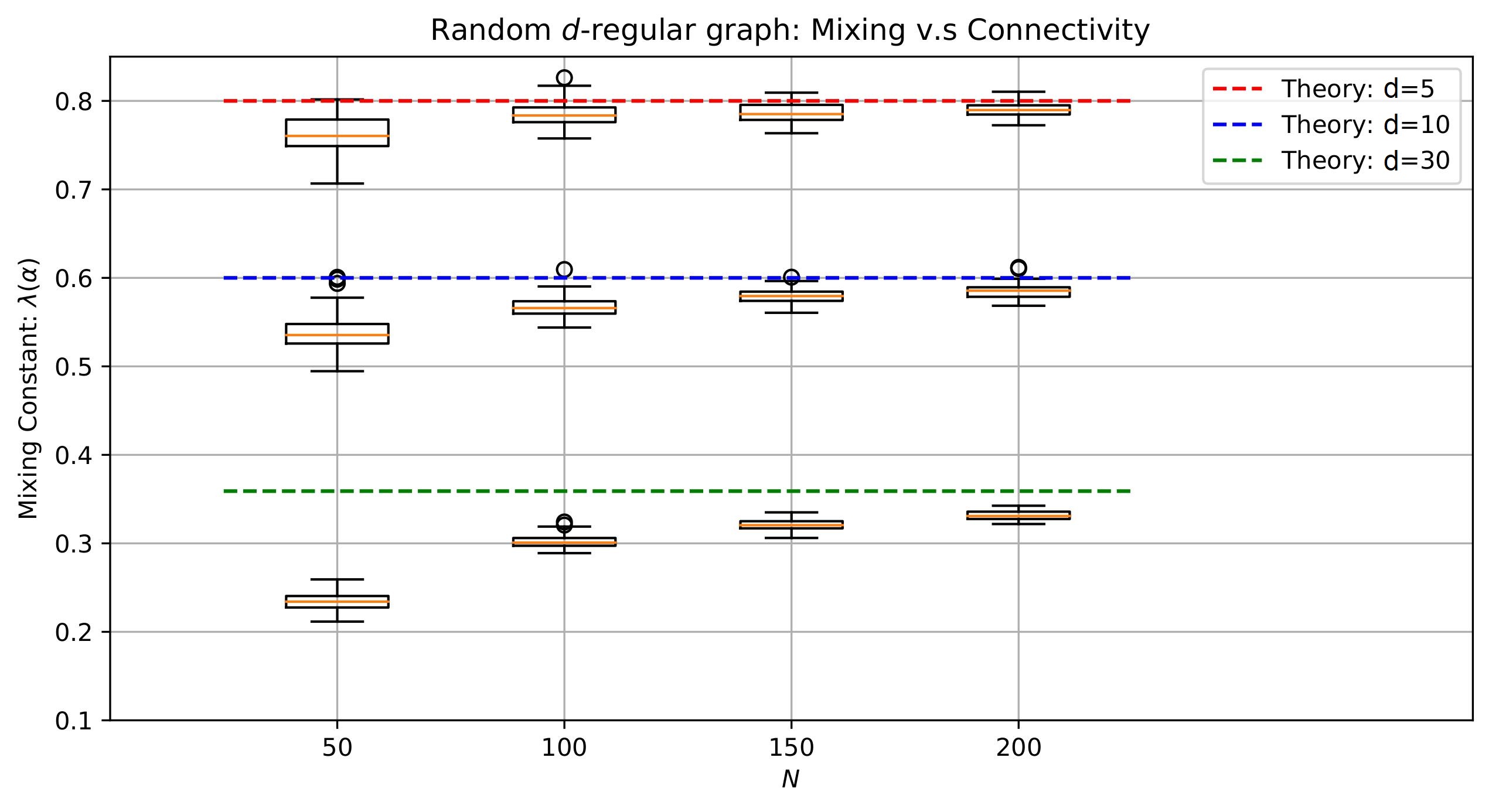}
\caption{Mixing constants $\lambda(\alpha)$ for several values of $d$ and $N$.}
\label{fig.spectral_gaps_1}
\end{figure}

We use the random graph generation algorithm of \cite{steger1999generating} as implemented in the \texttt{NetworkX} package of Python \citep{hagberg-2008-exploring} to generate random $\alpha$-matrices. This generates graphs $\G_N$ that are samples from the uniform distribution over all $d$-regular graphs with $N$ nodes. The $\alpha$-matrix is defined as the Markov transition matrix of the random walk on $\G_N$. 
We consider different values of $(d, N)$ and simulate 100 random $\alpha$ matrices for each pair. 
\cref{fig.spectral_gaps_1} shows the quality of the mixing constant $\lambda(\alpha)$ as a function of $d$ and $N$, as well as the limiting value as $N \to \infty$ as described in \cref{eq.alon.friedman}.

%
%
\subsection{Predictive distribution estimations}

\begin{figure}[!ht]
\centering
\includegraphics[width=0.7\textwidth]{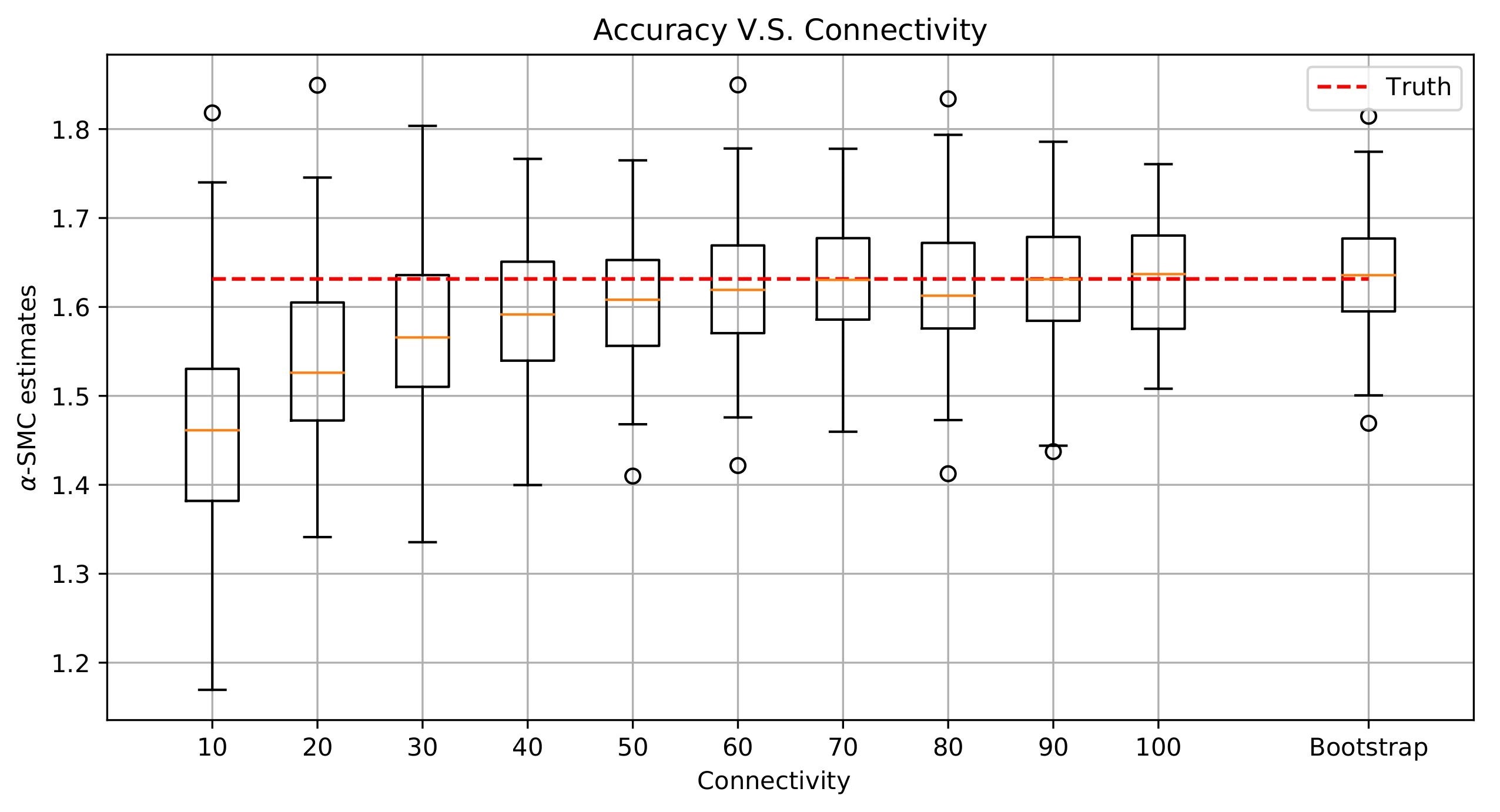} 
\includegraphics[width=0.7\textwidth]{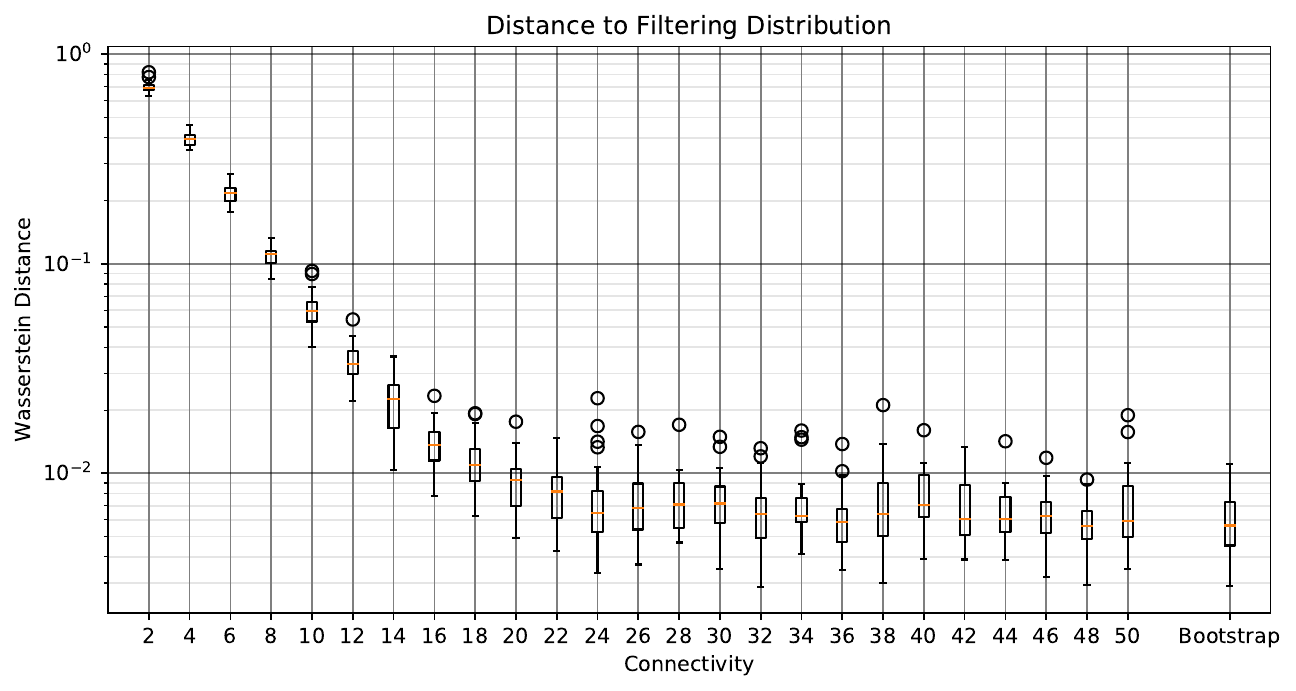}
\caption{The top plot displays the accuracy of the estimation of $\pi_{T=6}(\phi)$ with $\phi(x) = x^2$ and $N=10^4$ particles.
The bottom plot displays the Wasserstein distance between the estimated predictive distributions and the truth, which measures the accuracy of the estimated predictive distribution at time index $T=6$.}
\label{fig.wasserstein}
\end{figure}

Consider the state-space-model with initial distribution $\pi_0 = \N(0,1)$, dynamics $X_{t+1} = \beta X_{t} + \sqrt{1-\beta^2} \xi_t$ with $\beta=0.9$ and $\xi_t \sim \N(0,1)$. We consider the situation where the potential functions are all equal and given by $g_t(x) = 0.1 + 10 \times \indicator (|x-2|<0.1)$ for $t \geq 0$. We run several experiments with $N = 10^4$ particles for different values of the connectivity $d$. 
For each experiment, we randomly generate a $d$-regular graph as described in \cref{sec.randomalpha} and run the $\alpha$-sequential Monte Carlo algorithm using this. 
The top panel of \cref{fig.wasserstein} shows the performance of the $\alpha$-sequential Monte Carlo algorithm for the estimation of $\pi_{T=6}(\phi)$ for $\phi(x) = x^2$. For a connectivity $d=50$, the estimate from $\alpha$-sequential Monte Carlo is roughly as accurate as the bootstrap particle filter. The bottom panel of \cref{fig.wasserstein} shows the Wasserstein distance between the estimated predictive distributions and the true predictive distribution obtained by running an $\alpha$-sequential Monte Carlo algorithm for several values of the connectivity $d \geq 0$; the true predictive distribution is obtained by running the bootstrap particle filter with a large number of particles.

\subsection{Comparison with bootstrap particle filter} \label{sec.alphavsbootstrap}

We numerically investigate the effects of using sparse $d$-regular networks on the stability of the $\alpha$-sequential Monte Carlo algorithm. Three settings are considered.
\begin{itemize}
\item[(a)]
A local exchange scheme \citep{heine2017fluctuations}.

\item[(b)]
Generating an $\alpha$ matrix as described in \cref{sec.randomalpha} at the beginning of the algorithm and fixing it throughout. This is the setting considered in \cref{sec.non-asymptotic_analysis} and is referred to as `random $d$-regular (no permutation)' in this section.

\item[(c)]
Randomly permuting the matrix generated in (b) at time time step; this is the setting considered in \cref{sec:randomized_connections} and is referred to as ``random $d$-regular (with permutation)' in this section.

\end{itemize}
We consider a time-discretized version of the chaotic Lorenz 63 model \citep{lorenz1963deterministic}. The hidden chain $\{X_t\}_{t \geq 0}$ is three-dimensional with $X_t = (X_{t,1},X_{t,2},X_{t,3})$ and evolves as 
\begin{align*}
X_{t+\Delta t,1} 
& =
X_{t,1} + \Delta t \sigma (X_{t,2}-X_{t,1}) + \varepsilon_{t,1},
\\
X_{t+\Delta t,2} 
& = 
X_{t,2} + \Delta t \{ X_{t,1}(\rho - X_{t,3}) - X_{t,2} \} + \varepsilon_{t,2},
\\
X_{t+\Delta t,3}
& = 
X_{t,3} + \Delta t (X_{t,1} X_{t,2} - \beta X_{t,3}) + \varepsilon_{t,3},
\end{align*}
where $\Delta t = 10^{-3}$ is the time-discretization and $\varepsilon_t = (\varepsilon_{t,1}, \varepsilon_{t,2}, \varepsilon_{t,3})$ are independent and identically distributed as $\N(0, \Delta t \tau^2 I)$ for $\tau = 10^{-1}$. This model is known to be chaotic when $(\sigma,\rho,\beta) = (10,28,8/3)$, and this is the setting we choose. 
We collect observations $Y_t$ after every $\delta = 10 \Delta t$ units of time and assume that they are distributed as $Y_t \mid X_t \sim \N(X_t, \eta^2 I)$ for $\eta = 5 \times 10^{-1}$.

We generate $T=10^3$ observations from this model. The bootstrap particle filter with $10^6$ particles is used to calculate the ground truth. 
We compare the relative mean square errors of the estimate to the log-likelihood and predictive mean $\EE(X_T \mid Y_{0:(T-1)})$ for the three methods; this is the ratio of the mean square error of the estimate obtained by each method to the mean square error of the estimate obtained by the bootstrap particle filter with the same number of particles. We repeat the experiments $100$ times to obtain the mean square error.

\begin{figure}
\centering
\includegraphics[width=0.49\textwidth]{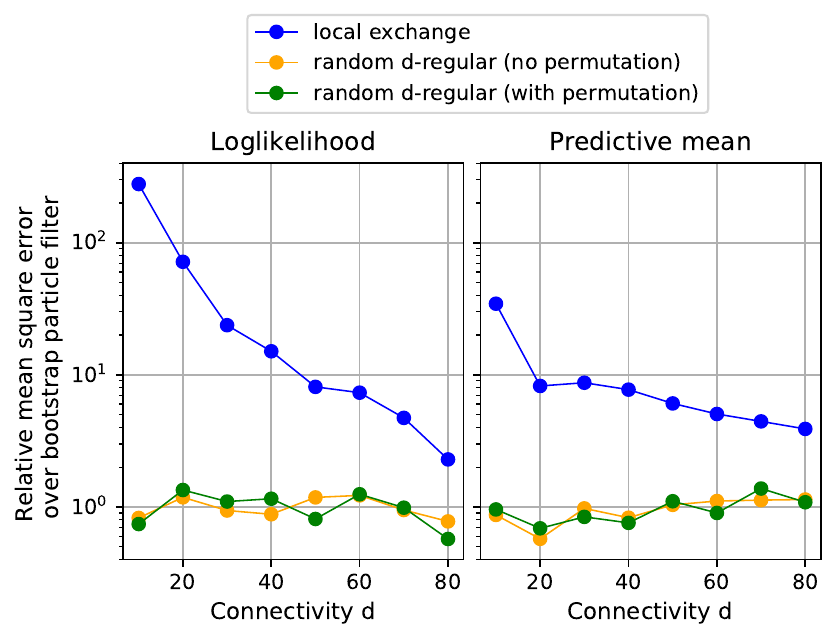}
\includegraphics[width=0.49\textwidth]{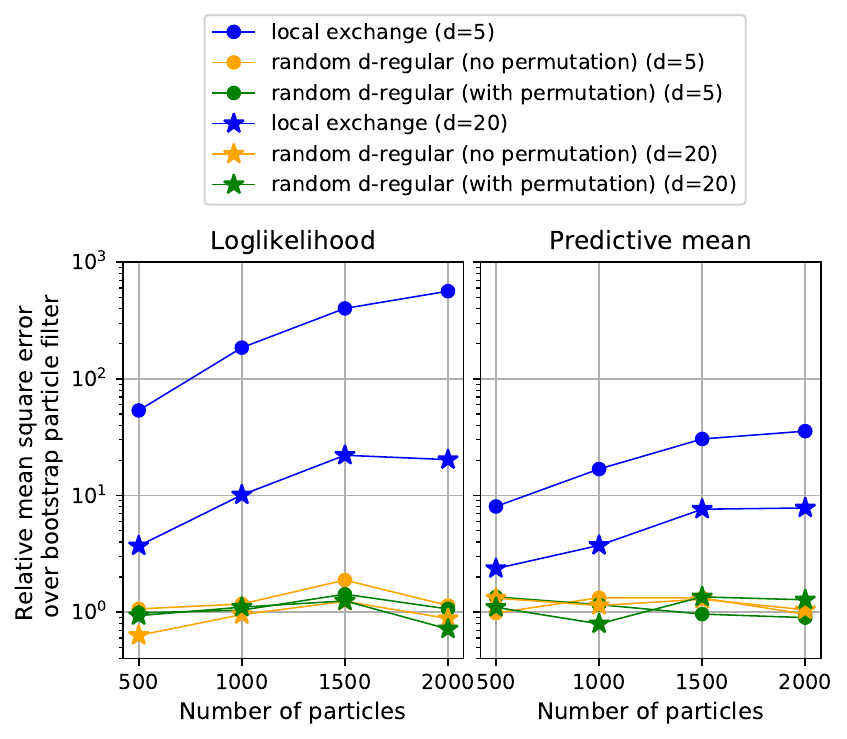} 
\caption{Relative performances with respect to the bootstrap particle filter; the left two plots are for $N=5\times10^3$.}
\label{expts.particles+connections}
\end{figure}

The two left plots of \cref{expts.particles+connections} display relative mean square errors for $N=5 \times 10^4$ as the degree $d$ of the graph increases.
As expected, the local exchange particle filter has a large error as compared to the bootstrap particle filter, which decreases as the degree increases. More interestingly, choosing a random $d$-regular graph has much lower error and is virtually indistinguishable from the bootstrap particle filter. This is true irrespective of whether we permute the nodes of the graph at every time, which is unsurprising as the permutation operation leaves the mixing constant unchanged. A random $5$-regular graph appears to perform extremely well. 

The two right plots of \cref{expts.particles+connections} display relative mean square errors as the network size (number of particles) $N$ increases. The performance of the local exchange particle filter deteriorates as $N$ increases, which is unsurprising since its mixing deteriorates. However, as predicted by the theory, the performance of a random $d$-regular graph remains stable as $N$ increases, whether or not the nodes are permuted at each time.

\begin{figure}
\centering
\includegraphics[width=0.9\textwidth]{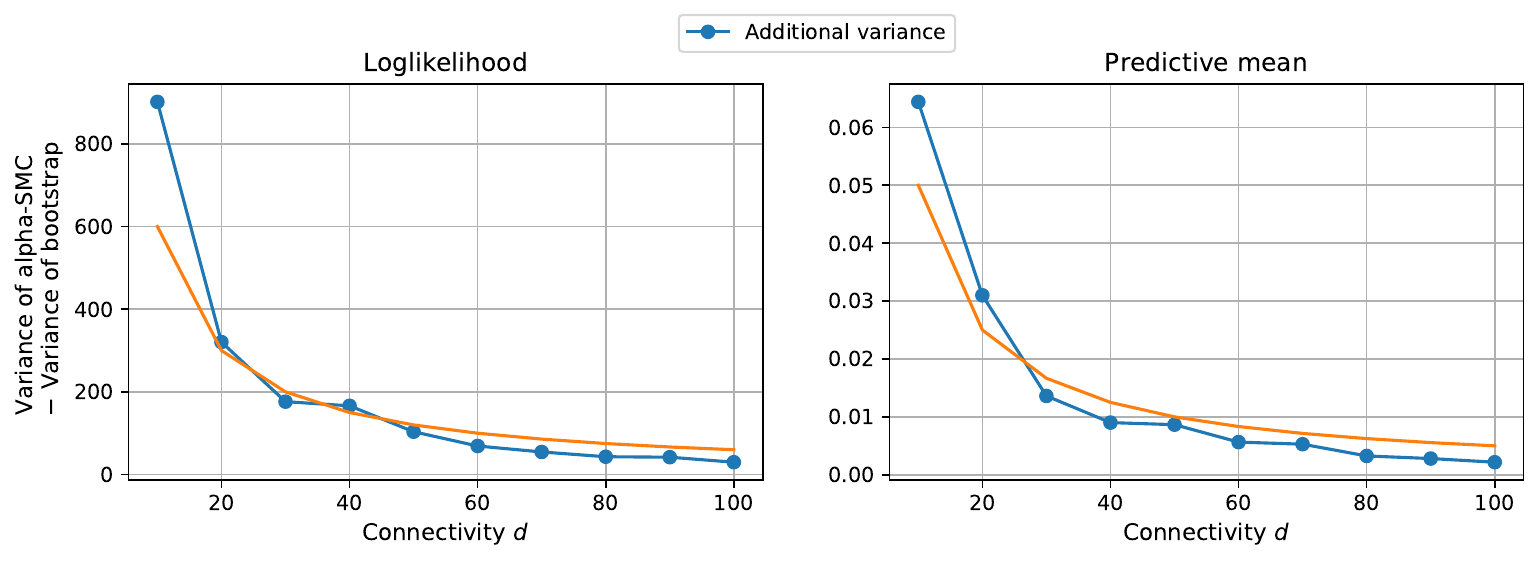}
\caption{Additional variance of the $\alpha$-sequential Monte Carlo algorithm as compared to the bootstrap particle filter for $N = 5 \times 10^3$. The orange line is proportional to $d^{-1}$.}
\label{expts.alphaSMC_v_bootstrap}
\end{figure}

Finally, we display in \cref{expts.alphaSMC_v_bootstrap} the additional variance of the $\alpha$-sequential Monte Carlo algorithm as compared to the bootstrap particle filter when using a random $d$-regular graph as the connectivity structure. As predicted by the theory, the additional variance is of order $\O(d^{-1})$.

\section{Conclusion} \label{sec.conclusion}

The bottleneck in parallelising particle filters is usually the resampling step since it typically involves interactions between all particles. Reducing these interactions can lead to more efficient algorithms, albeit sometimes at the expense of stability. 
Future directions can include relaxing the assumptions made in this article, considering adaptive sequential Monte Carlo \citep{fearnhead2013adaptive}, and considering high-dimensional target spaces \citep{beskos2014stability}. An interesting future direction would be to consider estimating the variance of the estimates obtained by the $\alpha$-sequential Monte Carlo algorithm along the lines of \cite{chan2013general,lee2018variance}.
From an applied perspective, it would be interesting to compare stable distributed implementations of sequential Monte Carlo with distributed Markov chain Monte Carlo.

\section*{Acknowledgment} 
The author acknowledges support from grant DMS-1638521 from SAMSI. The author would like to thank Alexandre H Thiery and Kari Heine for helpful discussions.

\appendix

\section{Proof of time-uniform stability} \label{proof:non-asymptotic_stability}

\begin{proof}[Proof of \cref{thm:non-asymptotic_stability}]

Recall that the sequence of probability measures $\{ \pi_t \}_{t \geq 0}$ defined by the state-space model \eqref{eq:normalised_target} of the main text satisfies $\pi_t = \forward_{s,t} \pi_s$,
where the operator $\forward_t$ is defined as
\begin{align*}
\forward_t \pi (\varphi) & = \frac{ \pi ( g_{t-1} \kernel_t \varphi ) }{ \pi (\potential_{t-1}) } \quad (t \geq 1).
\end{align*}

The stability properties of the operators $\{ \forward_t \}_{t \geq 1}$ are well-understood \citep{del2004feynman}. Under \cref{assumption_1} of the main text, there exist constants $D > 0$ and $\rho \in (0,1)$ such that for any two probability measures $\mu, \mu' \in \P(\XX)$, we have $\tnorm{ \forward_{s,t} \mu - \forward_{s,t} \mu' } \; \leq \; D \rho^{t - s} \tnorm{ \mu - \mu' }$.
%
The decomposition $(\widehat{\pi}^N_T - \pi_T) =( \widehat{\pi}^N_T - \forward_{0,T} \widehat{\pi}^N_0) + ( \forward_{0,T} \widehat{\pi}^N_0 - \forward_{0,T} \pi_0 )$ and the standard telescoping expansion $(\widehat{\pi}^N_T - \forward_{0,T} \widehat{\pi}^N_0) = \sum_{t=1}^T ( \forward_{t,T} \widehat{\pi}^N_t - \forward_{t,T} \forward_t \widehat{\pi}^N_{t-1} )$ yields that the discrepancy $\tnorm{ \widehat{\pi}^N_T - \pi_T }$ can be controlled as
\begin{align} \label{eq.error_decomposition}
\begin{aligned} 
\tnorm{ \widehat{\pi}^N_T - \pi_T }
&\leq
\sum_{t=1}^T \tnorm{ \forward_{t,T} \widehat{\pi}^N_t - \forward_{t,T} \forward_t \widehat{\pi}^N_{t-1} }
\; + \; 
\tnorm{ \forward_{0,T} \widehat{\pi}^N_0 - \forward_{0,T} \pi_0 } \\
&\leq 
\sum_{t=1}^T D \rho^{T-t} \tnorm{ \widehat{\pi}^N_t - \forward_t \widehat{\pi}^N_{t-1} }
\; + \; 
D \rho^T \tnorm{ \widehat{\pi}^N_0 - \pi_0 }.
\end{aligned}
\end{align}
Since $\widehat{\pi}^N_0 (\phi) = N^{-1} \sum_{i=1}^N \phi(X_{0,i})$ for independent and identically distributed samples $X_{0,i} \sim \pi_0$, it follows that $\tnorm{ \widehat{\pi}^N_0 - \pi_0 } \leq N^{-1/2}$. Consequently, since $\rho \in (0,1)$, for proving an upper bound of the type $\sup_{t \geq 0} \; \tnorm{\widehat{\pi}^N_t - \pi_t} \; \leq \; \textrm{Cst} \times N^{-1/2}$, it only remains to prove that the quantities $\tnorm{ \widehat{\pi}^N_t - \forward_t \widehat{\pi}^N_{t-1} }$ can be uniformly bounded in time by a constant multiple of $N^{-1/2}$. 

For a test function $\varphi \in \B(\XX)$, we have $\widehat{\pi}^N_t (\varphi) - \forward_t \widehat{\pi}^N_{t-1}(\varphi) = (\widetilde{A}- A)/B$ for \\ $A = \sum_{i=1}^N W_{t-1}^i \potential_{t-1}(X_{t-1}^{i}) \kernel_t \varphi (X_{t-1}^{i})$, $\widetilde{A} = \sum_{i=1}^N W_{t}^i \varphi(X_{t}^{i})$, and $B = \sum_{i=1}^N W_{t-1}^i \potential_{t-1}(X_{t-1}^{i})$.
We have $\EE_{t-1}(\widetilde{A}) = A$, and the quantities $A$, $B$ and $W^i_t$ are all $\F_{t-1}$-measurable. It follows that
\begin{align} \label{eq.link.stability.ess}
\begin{aligned}
& ~~ 
\EE_{t-1} \sqBK{ \curBK{ \widehat{\pi}^N_t (\varphi) - \forward_t \widehat{\pi}^N_{t-1}(\varphi) }^2 }
=
B^{-2} \Var_{t-1} ( \widetilde{A} )
\\
& =
B^{-2} \sum_{i=1}^N ( W_{t}^i )^2 \Var_{t-1} \curBK{ \varphi(X_{t}^{i}) } 
\leq
B^{-2} \sum_{i=1}^N (W_{t}^i)^2 
= \norm{ \overline{W}_t }^2 = \mathsf{E}^N_t.
\end{aligned}
\end{align}
In the last line of \cref{eq.link.stability.ess}, we have used the fact that $B = \sum_{i=1}^N W_{t}^i$. We have also introduced the quantity $\mathsf{E}^N_t = \| \overline{W}_t \|^2$; this is a measure of the effective sample size \citep{whiteley2016role}. In summary, we have thus established that 
\begin{align} \label{eq.ess_norm}
\tnorm{ \widehat{\pi}_t^N - \forward_t \widehat{\pi}_{t-1}^N }^2 & \leq \EE \left ( \mathsf{E}_t^N \right ).
\end{align}

As recognised in \cite{whiteley2016role}, \cref{eq.link.stability.ess} shows that controlling the behaviour of $\mathsf{E}_t^N$ is crucial to studying the stability properties of the $\alpha$-sequential Monte Carlo algorithm. 
For proving a bound of the type given by $\sup_{t \geq 0} \; \tnorm{\widehat{\pi}^N_t - \pi_t} \; \leq \; \textrm{Cst} \times N^{-1/2}$, \cref{eq.ess_norm} reveals that it suffices to have the uniform-in-time bound $\EE (\mathsf{E}_t^N) \leq \mathrm{Cst} / N$. Recalling that $\kappa_g^{-1} \leq g_{t-1}(x) \leq \kappa_g$ by \cref{assumption_1} of the main text, the bound \eqref{eq.fundamental.bound} yields that
\begin{align} \label{eq.recursion.E.bound}
\begin{aligned}
\mathsf{E}_t^N
&=
\frac{\sum_{i=1}^N \{ \sum_{j=1}^N \alpha^{ij} W_{t-1}^{j} g_{t-1}(X_{t-1}^j) \}^2}{ \{ \sum_{i=1}^N \sum_{j=1}^N \alpha^{ij} W_{t-1}^{j} g_{t-1}(X_{t-1}^j) \}^2}
\leq
\kappa_g^4 \sum_{i=1}^N \left ( \sum_{j=1}^N \alpha^{ij} \overline{W}_{t-1,j} \right )^2
=
\kappa_g^4 \norm{\alpha \overline{W}_{t-1}}^2 \\
& \leq 
\kappa_g^4 \curBK{ \frac{1-\lambda^2(\alpha)}{N} + \lambda^2(\alpha) \norm{\overline{W}_{t-1}}^2 } 
= 
\kappa_g^4 \curBK{ \frac{1-\lambda^2(\alpha)}{N} + \lambda^2(\alpha) \mathsf{E}_{t-1}^N}.
\end{aligned}
\end{align}
If the constant $\lambda(\alpha)$ defined in \cref{eq.constant.lambda} of the main text satisfies $\lambda(\alpha) < 1/\kappa_g^{2}$, iterating the bound \eqref{eq.recursion.E.bound} directly yields that
\begin{align} \label{eq.E.bound}
\tnorm{ \widehat{\pi}_t^N - \forward_t \widehat{\pi}_{t-1}^N }^2 & \leq \EE \BK{ \mathsf{E}_t^N }
\leq \frac{\kappa_g^4 \{1 - \lambda^2(\alpha)\}}{1 - \kappa_g^4 \lambda^2(\alpha)} \frac{1}{N}
\quad \textrm{for all} ~ t \geq 0.
\end{align}
Combining \cref{eq.error_decomposition} and \cref{eq.E.bound}, the theorem is proved.
\end{proof}

\section{Proofs for randomized connections}

\subsection{Setup}
\label{sec.random.connect.setup}

The following proposition is useful in studying the asymptotic behaviour of the $\alpha$-sequential Monte Carlo algorithm. 
\begin{proposition}[Basic properties] \label{prop.basic}
The following are true, where the expectations are with respect to the distribution $\alphadist_N$ on $\bistochm_N$.
\begin{enumerate}[(a)]

\item \label{prop.a}
$\EE\{(\alpha^{ii})^2\}$ does not depend on $i$.

\item \label{prop.b}
For $i \neq j$, $\EE(\alpha^{ij})$ and $\EE\{(\alpha^{ij})^2\}$ do not depend on $(i,j)$. Further, there exists $0\leq c_1 < \infty$ such that $\EE(\alpha^{ij}) \leq c_1 N^{-1}$ and $\EE\{(\alpha^{ij})^2\} \leq c_1 N^{-1}$ for $N$ large.

\item \label{prop.c}
For $i \neq j$, $\EE(\alpha^{ii} \alpha^{ij})$ does not depend on $(i,j)$, Further, there exists $0 \leq c_2 < \infty$ such that $\EE(\alpha^{ii}\alpha^{ij} ) \leq c_2 N^{-1}$ for $N$ large.
 
\end{enumerate}
\end{proposition}

For example, for $i \neq j$, we have $\EE (\alpha^{ij}) = N^{-1}$ and $\EE \{ (\alpha^{ij})^2 \} = N^{-2}$ for the bootstrap particle filter, and we have $\EE(\alpha^{ij}) = \EE \{ (\alpha^{ij})^2 \} = 0$ for sequential importance sampling. 

Let the permutation corresponding to a permutation matrix $P$ be $\sigma : \{1,\dots,N\} \mapsto \{1,\dots,N\}$. Then $(P \alpha P^{-1})^{ij} = \alpha^{\sigma(i) \sigma^{-1}(j)}$.
By \cref{eq.permute.property}, this implies $\alpha^{ij} \eqd \alpha^{\sigma(i) \sigma^{-1}(j)}$ for all $1 \leq i,j \leq N$ and all permutations $\sigma$, where we have used the notation $\cdot \eqd \cdot$ to denote that the left and right hand side have the same distribution.

\begin{proof}[Proof of \cref{prop.basic}]

Part \ref{prop.a} follows since for any $i \neq i'$, there exists a permutation $\sigma$ such that $\sigma(i) = i'$ and $\sigma(i')=i$, which implies that $\alpha^{ii} \eqd \alpha^{\sigma(i) \sigma^{-1}(i)} = \alpha^{i' i'}$.

To see part \ref{prop.b}, consider $i \neq j \neq k$. There exists a permutation $\sigma$ such that $\sigma(i) = i$ and $\sigma(k) = j$, which implies that $\alpha^{ij} \eqd \alpha^{\sigma(i) \sigma^{-1}(j)} = \alpha^{ik}$. This implies that $\EE(\alpha^{ij}) = \EE(\alpha^{ik})$ for all $i \neq j \neq k$. Similarly, we have $\EE(\alpha^{i j}) = \EE(\alpha^{i' j})$ for all $i \neq i' \neq j$. The previous two statements imply that $\EE(\alpha^{ij})$ does not depend on $(i,j)$ for $i \neq j$.
Since $\sum_{j=1}^N \alpha^{ij} = 1$ and $\alpha^{ij}\leq1$, the results follow. 

To see part \ref{prop.c}, consider $i \neq j \neq k$. There exists a permutation $\sigma$ such that $\sigma(i)=i$ and $\sigma(k)=j$, and therefore $(\alpha^{ii}, \alpha^{ij}) \eqd (\alpha^{\sigma(i) \sigma^{-1}(i)}, \alpha^{\sigma(i) \sigma^{-1}(j)}) = (\alpha^{ii}, \alpha^{ik})$ for $j \neq k$. Thus $\EE(\alpha^{ii} \alpha^{ij})$ does not depend on $j$ for $j \neq i$. Similarly, for $i \neq i' \neq j$, there exists a permutation $\sigma$ such that $\sigma(i)=i'$, $\sigma(i')=i$, and $\sigma(j)=j$, which implies $(\alpha^{ii}, \alpha^{ij}) \eqd (\alpha^{\sigma(i)\sigma^{-1}(i)},\alpha^{\sigma(i) \sigma^{-1}(j)}) = (\alpha^{i'i'}, \alpha^{i'j})$. Thus $\EE(\alpha^{ii} \alpha^{ij})$ does not depend on $i$ for $i \neq j$.
The inequality follows from part \ref{prop.b} as $\alpha^{ij} \leq 1$.
\end{proof}

By the assumption that for any time index $t \geq 0$ and $x \in \XX$, we have $0< g_t(x) \leq \kappa_g$, it follows that for any time index $t \geq 0$ and particle index $1 \leq i \leq N$, we have $0 < W^i_t \leq \kappa_g^t$,
so that, for a test function $\phi \in \mathcal{B}(\XX)$, the random variables $\gammahat^N_t(\phi)$ and $\muhat^N_t(\phi)$ are almost surely bounded: $\|\gammahat^N_t(\phi)\|_{\infty} \leq \kappa_g^t$ and $\|\muhat^N_t(\phi)\|_{\infty} \leq \kappa_g^{2t}$, where the infinity norm of a random variable $X$ is defined as $\|X\|_\infty = \inf \{B > 0 : \PP(|X| \leq B) = 1 \}$.

For the bootstrap particle filter, it is standard that as $N \to \infty$, the sequence of particle approximations $\pihat^N_t$ is consistent \citep{del1996non}. Since in that case, all the weights $W^i_t$ are equal and converge in probability to $Z_t$, it follows that
\begin{equation} \label{eq.bootstrap.pf.mu.flow}
\muhat_t^N(\varphi)
=
\frac1N \sum_{i=1}^N (W^i_t)^2 \phi(X^i_t) 
\plim
Z_t^2 \pi_t(\phi),
\end{equation}
where $\plim$ denotes convergence in probability. Equation \eqref{eq.flow.mu.basic} of the main text can be seen as a generalisation of \cref{eq.bootstrap.pf.mu.flow}.

For the bootstrap particle filter, \cref{eq.flow.mu.basic} of the main text implies $\mu_t(\phi) = Z_t^2 \pi_t(\phi)$, which in turn implies \cref{eq.bootstrap.pf.mu.flow}. Similarly, for importance sampling, \cref{eq.flow.mu.basic} of the main text implies $\mu_t(\phi)=\mu_{t-1}(\forwardQQ_t \phi)$, which is as expected.
Equation \eqref{eq.flow.mu.basic} is therefore a generalization of the bootstrap particle filter and importance sampling, which represent two extreme communication structures (fully connected and not connected at all, respectively).

\subsection{Consistency} 
\label{proof:consistency}

\begin{proof}[Proof of \cref{thm.consistency}]

Since $\widehat{\pi}^N_t(\phi) = \gammahat^N_t(\phi) / \gammahat^N_t(1)$, it suffices to prove that $\gammahat^N_t(\phi) \plim \gamma_t(\phi)$ and $\muhat^N_t(\phi) \plim \mu_t(\phi)$. We work by induction, the initial case $t=0$ following directly from the weak law of large numbers.\\

\noindent
{\bf (I) Convergence of $\gammahat^N_t$.}\\ 
This follows from Theorem 1 of \cite{whiteley2016role}, but we include a proof for the sake of completeness. 
Since $\|\gammahat^N_t(\phi)\|_{\infty} \leq \kappa_g^t$ and $\EE\{\gammahat^N_t(\phi)\} = \gamma_t(\phi)$, for proving $\gammahat^N_t(\phi) \plim \gamma_t(\phi)$, it suffices to prove that the variance of $\gammahat^N_t(\phi)$ converges to zero. 

\begin{itemize}

\item
The variance of $\EE_{t-1}\{\gammahat^N_t(\phi)\} = \gammahat^N_t(\forwardQ_t \phi)$ converges to zero since the random variables $\gammahat^N_t(\forwardQ_t \phi)$ are bounded by $\kappa_g^t$ and, by induction, converge in probability to $\gamma_t(\forwardQ_t \phi)$ as $N \to \infty$.

\item 
The expectation of $\Var_{t-1}\{\gammahat^N_t(\phi)\}$ also converges to zero. Indeed, since the random variables $\{W^i_t \, \phi(X^i_t)\}_{i=1}^N$ are independent conditionally upon $\F_{t-1}$ and upper bounded in absolute value by $\kappa_g^t$, we have that, almost surely,
\begin{align*}
\Var_{t-1} \curBK{\frac{1}{N} \sum_{i=1}^N W^i_t \, \phi(X^i_t)} \leq \frac{\kappa_g^{2t}}{N} \to 0.
\end{align*}
\end{itemize}
Since $\Var\{ \gammahat^N_t(\phi) \}$ can be decomposed as the sum of the expectation of $\Var_{t-1}\{\gammahat^N_t(\phi)\}$ and the variance of $\EE_{t-1}\{\gammahat^N_t(\phi)\}$, this concludes the proof that $\gammahat^N_t(\phi) \plim \gamma_t(\phi)$.\\

\noindent
{\bf (II) Convergence of $\muhat^N_t$.}\\ 
We proceed in two steps. We first show that $\EE\{\muhat^N_t(\phi)\} \to \mu_t(\phi)$, and then prove that $\Var\{\muhat^N_t(\phi)\}$ converges to zero. Since $\|\muhat^N_t(\phi)\|_{\infty} \leq \kappa_2^{2t}$, this is enough to obtain that $\muhat^N_t(\phi) \plim \mu_t(\phi)$.

To begin with,
\begin{align}
& \quad 
\EE_{t-1} \curBK{\muhat^N_t(\phi)} 
= 
\frac{1}{N} \sum_{i=1}^N \EE_{t-1} \curBK{(W^i_t)^2 \, \phi(X^i_t)}
\nonumber \\
& =
\frac1N \sum_{i=1}^N \EE_{t-1} \sqBK{ \curBK{\sum_{j=1}^N \alpha^{ij}_{t-1} \, W^{j}_{t-1} \, g_{t-1}(X^{j}_{t-1})} \times \curBK{\sum_{k=1}^N \alpha^{ik}_{t-1} \, W^{k}_{t-1} \, g_{t-1}(X^{k}_{t-1}) \, K_t \phi(X^{k}_{t-1})} }
\nonumber \\
& = 
\frac1N \sum_{i=1}^N \EE_{t-1} \curBK{\sum_{j=1}^N (\alpha_{t-1}^{ij})^2 (W_{t-1}^j)^2 g^2_{t-1}(X_{t-1}^{j}) K_t \phi(X_{t-1}^j)}
\label{eq.1}
\\
& \qquad 
+ \frac1N \sum_{i=1}^N \EE_{t-1} \curBK{\sum_{j=1}^N \sum_{\substack{k=1 \\ k \neq j}}^N \alpha_{t-1}^{ij} \alpha_{t-1}^{ik} W^{k}_{t-1} \, W^{j}_{t-1} \, g_{t-1}(X^{k}_{t-1}) \, g_{t-1}(X^{j}_{t-1}) \, K_t \phi(X^{k}_{t-1})}.
\label{eq.2}
\end{align}

Expression \eqref{eq.1} is 
\begin{align*}
& \quad 
\frac1N \sum_{i=1}^N \EE \curBK{(\alpha^{ii})^2} (W_{t-1}^i)^2 g^2_{t-1}(X_{t-1}^{i}) K_t \phi(X_{t-1}^i)
\\
& \quad 
+ 
\frac1N \sum_{i=1}^N \sum_{\substack{j=1 \\ j \neq i}}^N \EE \{(\alpha^{ij})^2\} (W_{t-1}^j)^2 g^2_{t-1}(X_{t-1}^{j}) K_t \phi(X_{t-1}^j)
\\
& = 
\frac{\EE \{ (\alpha^{11})^2 \}}{N} \sum_{i=1}^N (W_{t-1}^i)^2 g^2_{t-1}(X_{t-1}^{i}) K_t \phi(X_{t-1}^i)
\\
& \quad + 
\frac{\EE \{(\alpha^{12})^2\}}{N} \sum_{i=1}^N \sum_{\substack{j=1 \\ j \neq i}}^N  (W_{t-1}^j)^2 g^2_{t-1}(X_{t-1}^j) K_t \phi(X_{t-1}^j)
\\
& = 
\EE \{ (\alpha^{11})^2 \} \times \muhat^N_{t-1}(\forwardQQ_t \phi)
\\
& \quad 
+ 
\frac{\EE\{(\alpha^{12})^2\}}{N} \sum_{i=1}^N \curBK{\sum_{j=1}^N (W_{t-1}^j)^2 g^2_{t-1}(X_{t-1}^j) K_t \phi(X_t^j) - (W_{t-1}^i)^2 g^2_{t-1}(X_{t-1}^i) K_t \phi(X_{t-1}^i)} 
\\
& = 
\EE \{ (\alpha^{11})^2 \} \times \muhat^N_{t-1}(\forwardQQ_t \phi)
\\
& \quad +
N \EE\{(\alpha^{12})^2\} \left \{ \frac1N \sum_{j=1}^N (W_{t-1}^j)^2 g^2_{t-1}(X_{t-1}^j) K_t \phi(X_t^j) + \O \left ( \frac1N \right ) \right \}
\\
& = 
\muhat^N_{t-1}(\forwardQQ_t \phi) \times \sqBK{\EE \{ (\alpha^{11})^2 \} + N \EE\{(\alpha^{12})^2\}}
+
\O \left ( \frac1N \right )
\\
& \plim 
\mu_{t-1}(\forwardQQ_t \phi) \times \lim_{N \to \infty} \sqBK{\EE \{(\alpha^{11})^2\} + N \EE\{(\alpha^{12})^2\}},
\end{align*}
where the first equality is by \cref{prop.basic}\ref{prop.a} and the first part of \cref{prop.basic}\ref{prop.b}, 
the second and fourth equalities are because $\muhat^N_{t-1}(\forwardQQ_t \phi) = (1/N) \sum_{i=1}^N (W_{t-1}^i)^2 g^2_{t-1}(X_{t-1}^{i}) K_t \phi(X_{t-1}^i)$, the third equality is by the second part of \cref{prop.basic}\ref{prop.b}, and
and the limit is by the consistency of $\muhat^N_{t-1}$. 

For the sake of convenience, define 
\begin{align*}
b_t^{jk} 
& =
W^{k}_{t-1} \, W^{j}_{t-1} \, g_{t-1}(X^{k}_{t-1}) \, g_{t-1}(X^{j}_{t-1}) \, K_t \phi(X^{k}_{t-1}),
\\
a_t^{ijk} 
& =
\EE(\alpha_{t-1}^{ij} \alpha_{t-1}^{ik}) W^{k}_{t-1} \, W^{j}_{t-1} \, g_{t-1}(X^{k}_{t-1}) \, g_{t-1}(X^{j}_{t-1}) \, K_t \phi(X^{k}_{t-1}) 
= 
\EE(\alpha_{t-1}^{ij} \alpha_{t-1}^{ik}) \, b_t^{jk}.
\end{align*}
Expression \eqref{eq.2} can be written as
\begin{align*}
& \quad 
\frac1N \sum_{i=1}^N \sum_{j=1}^N \sum_{\substack{k=1 \\ k \neq j}}^N a_t^{ijk}
=
\frac1N \sum_{i=1}^N \left ( \sum_{\substack{k=1 \\ k \neq i}}^N a_t^{iik} + \sum_{\substack{j=1 \\ j \neq i}}^N \sum_{\substack{k=1 \\ k \neq j}}^N a_t^{ijk} \right )
\\
& = 
\frac1N \sum_{i=1}^N \left \{ \sum_{\substack{k=1 \\ k \neq i}}^N a_t^{iik} + \sum_{\substack{j=1 \\ j \neq i}}^N \left ( \sum_{\substack{k=1 \\ k \neq j,i}}^N a_t^{ijk} + a_t^{iji} \right )  \right \}
\\ 
& =
\frac1N \sum_{i=1}^N \sum_{\substack{k=1 \\ k \neq i}}^N \EE(\alpha_{t-1}^{ii} \alpha_{t-1}^{ik}) \, b_t^{ik}
+ 
\frac1N \sum_{i=1}^N \sum_{\substack{j=1 \\ j \neq i}}^N \sum_{\substack{k=1 \\ k \neq j,i}}^N \EE(\alpha_{t-1}^{ij} \alpha_{t-1}^{ik}) \, b_t^{jk}
\\ 
& \quad + 
\frac1N \sum_{i=1}^N \sum_{\substack{j=1 \\ j \neq i}}^N \EE(\alpha_{t-1}^{ij} \alpha_{t-1}^{ii}) \, b_t^{ji}
\\
& =
\frac{\EE(\alpha^{11} \alpha^{12})}{N} \sum_{i=1}^N \sum_{\substack{k=1 \\ k \neq i}}^N b_t^{ik} 
+ 
\EE(\alpha^{12} \alpha^{13}) \frac1N \sum_{i=1}^N \sum_{\substack{j=1 \\ j \neq i}}^N \sum_{\substack{k=1 \\ k \neq j,i}}^N b_t^{jk}
+ 
\frac{\EE(\alpha^{12}\alpha^{11})}{N} \sum_{i=1}^N \sum_{\substack{j=1 \\ j \neq i}}^N b_t^{ji}
\\
& = 
\frac{\EE(\alpha^{11}\alpha^{12})}{N} \sum_{i=1}^N \left ( \sum_{k=1}^N b_t^{ik} - b_t^{ii} + \sum_{j=1}^N b_t^{ji} - b_t^{ii} \right )
\\
& \quad + 
\frac{\EE(\alpha^{12} \alpha^{13})}{N} \sum_{i=1}^N \left ( \sum_{j=1}^N  \sum_{\substack{k=1 \\ k \neq j,i}}^N b_t^{jk} - \sum_{\substack{k=1 \\ k \neq i}}^N b_t^{ik} \right )
\\
& = 
\frac{\EE(\alpha^{11}\alpha^{12})}{N} \sum_{i=1}^N \left ( \sum_{k=1}^N b_t^{ik} + \sum_{j=1}^N b_t^{ji} - 2b_t^{ii} \right )
\\
& \quad + 
\frac{\EE(\alpha^{12} \alpha^{13})}{N} \sum_{i=1}^N \left \{ \sum_{j=1}^N \left ( \sum_{k=1}^N b_t^{jk} - b_t^{jj} - b_t^{ji} \right ) - \left ( \sum_{k=1}^N b_t^{ik} - b_t^{ii}\right )\right \}
\\
& = 
2N \EE(\alpha^{11}\alpha^{12}) \, \gammahat^N_{t-1}(g_{t-1}) \, \gammahat^N_{t-1} (g_{t-1} K_t \phi)
\\
& \quad + 
N^2 \EE(\alpha^{12}\alpha^{13}) \, \gammahat^N_{t-1}(g_{t-1}) \, \gammahat^N_{t-1} (g_{t-1} K_t \phi)
+ 
\O \left ( \frac1N \right )
\\
& \plim 
\gamma_{t-1} (g_{t-1} K_t \phi) \gamma_{t-1}(g_{t-1}) \times \lim_{N \to \infty} \curBK{ N^2 \EE(\alpha_{t-1}^{12} \alpha_{t-1}^{13}) + 2 N \EE(\alpha_{t-1}^{11}\alpha_{t-1}^{12} )}
\\
& =
Z_t^2 \pi_t(\phi) \times \lim_{N \to \infty} \curBK{ N^2 \EE(\alpha_{t-1}^{12} \alpha_{t-1}^{13}) + 2 N \EE(\alpha_{t-1}^{11}\alpha_{t-1}^{12} )},
\end{align*}
where the third equality uses \cref{prop.basic}\ref{prop.c}, the fourth equality uses \cref{ass.cross.terms}, and the sixth equality uses the fact all the relevant quantities are bounded.

Putting together what we have obtained so far, we get
\begin{align*}
\EE_{t-1}\{\muhat^N_t(\phi)\} 
& \plim
\mu_{t-1}(\forwardQQ_t \phi) 
\times
\lim_{N \to \infty} \sqBK{\EE \curBK{(\alpha^{11})^2} + N \EE\{(\alpha^{12})^2\}}
\\
& \quad +
Z_t^2 \pi_t(\phi) 
\times 
\lim_{N \to \infty} \curBK{ N^2 \EE(\alpha^{12} \alpha^{13}) + 2 N \EE(\alpha^{11}\alpha^{12} )}
\\
& = \mu_t(\phi).
\end{align*}

It remains to prove that $\Var\{\muhat^N_t(\phi)\} = \Var[ \EE_{t-1}\{\muhat^N_t(\phi)\}] + \EE[\Var_{t-1}\{\muhat^N_t(\phi)\}] \to 0$. We now prove that each term converges to zero. From the proof of $\EE_{t-1}\{\muhat^N_t(\phi)\} \to \mu_t(\phi)$, we have 
\begin{align} \label{eq.conditoined.expect.mu}
\begin{aligned}
\EE_{t-1}\{\muhat^N_t(\phi)\}
& =
\muhat^N_{t-1}(\forwardQQ_t \phi) 
\sqBK{\EE \curBK{(\alpha_{t-1}^{11})^2} + N \EE\{(\alpha^{12})^2\} }
\\
& \quad + 
N^2 \EE(\alpha_{t-1}^{12} \alpha_{t-1}^{13}) \curBK{\gammahat^N_{t-1}(g_{t-1}) \, \gammahat^N_{t-1} (g_{t-1} K_t \phi) + \O(N^{-1})}
\\
& \quad +
N \EE(\alpha_{t-1}^{11}\alpha_{t-1}^{12} )\curBK{\gammahat^N_{t-1}(g_{t-1}) \, \gammahat^N_{t-1} (g_{t-1} K_t \phi) + \O(N^{-1})}.
\end{aligned}
\end{align}

\begin{itemize}
\item 
By \cref{prop.basic} and \cref{ass.cross.terms}, all terms on the right-hand side of \cref{eq.conditoined.expect.mu} are upper bounded by a universal constant and, by induction, converge in probability to a constant. Therefore the variance of $\EE_{t-1}\{\muhat^N_{t-1}(\phi)\}$ converges to zero.
\item The expectation of $\Var_{t-1}\{\muhat^N_t(\phi)\}$ also converges to zero. Indeed, since the random variables $\{ (W^i_t)^2 \, \phi(X^i_t)\}_{i=1}^N$ are independent conditionally upon $\F^N_{t-1}$ and upper bounded in absolute value by $\kappa_g^{2t}$, we have that, almost surely,
\begin{align*}
\Var_t\curBK{\frac{1}{N} \sum_{i=1}^N (W^i_N)^2 \, \phi(X^i_N)} \leq \frac{\kappa_g^{4t}}{N} \to 0.
\end{align*}
\end{itemize}
This concludes the proof.
\end{proof}

\subsection{Central limit theorem}
\label{proof:CLT}

\begin{proof}[Proof of \cref{thm.central limit theorem}]
Notice that $\{\widehat{\pi}^N_t(\phi) - \pi_t(\phi)\} = \{\gammahat^N_t(1)\}^{-1} \, [\gammahat^N_t \{\phi - \pi_t(\phi) \}]$ and $\gammahat^N_t(1) \plim Z_t$. Consequently, the recursive formula for the asymptotic variance $\VV^\pi_t(\phi)$ readily follows from the one describing $\VV^\gamma_t(\phi)$. We thus concentrate on proving the recursive formula for $\VV^\gamma_t(\phi)$. We proceed by induction and use a standard Fourier-theoretic approach. The initial case $t=0$ follows directly from the standard central limit theorem for independent and identically distributed random variables. We need to prove that for any $\xi \in \RR$ and $S^N_t = N^{1/2} \{\gammahat^N_t(\phi) - \gamma_t(\phi)\}$, we have $\EE\{\exp\BK{\iu \xi S^{N}_t}\} \plim \exp\{-\VV^\gamma_t(\phi) \xi^2 / 2\}$, where $\iu$ denotes the imaginary unit. 
We have
\begin{align*}
S^N_t 
& = 
N^{1/2} \sqBK{\gammahat^N_t(\phi) - \EE_{t-1}\{\gammahat^N_t(\phi)\}}
+
N^{1/2} \sqBK{\EE_{t-1}\{\gammahat^N_t(\phi)\} - \gamma_t(\phi)}
\\
& =
A^N_t(\phi) + B^N_t(\phi).
\end{align*}
Further, $B^N_t(\phi) = N^{1/2} \{ \gammahat^N_{t-1}(\forwardQ \phi) - \gamma_{t-1}(\forwardQ \phi)\}$, so the induction hypothesis yields that $\EE[\exp\{\iu \xi B^N_t(\phi)\}] \plim \exp\{- \VV^\gamma_{t-1}(\forwardQ_t \phi) \xi^2/2\}$. To conclude, it suffices to show that
\begin{align} \label{eq.crux.central limit theorem}
\EE_{t-1} \sqBK{ \exp\curBK{\iu \xi A^N_t(\phi)} } 
\plim 
\exp\sqBK{-\curBK{\mu_t(\phi^2) - Z^2_t \pi_{t}(\phi)^2} \xi^2 / 2},
\end{align}
since it then follows (by Slutsky's theorem) that 
\begin{align*}
& \quad 
\EE\curBK{\exp(\iu \xi S_t^N)} = \EE\sqBK{\EE_{t-1} \sqBK{ \exp\curBK{\iu \xi A^N_t(\phi)} } \times \exp\curBK{\iu \xi B^N_t(\phi)} }
\\
& \implies 
\exp \sqBK{ - \curBK{\VV^\gamma_{t-1}(\forwardQ_t \phi) + \mu_t(\phi^2) - Z^2_t \pi_{t}(\phi)^2} \frac{\xi^2}{2}} 
= 
\exp\curBK{-\VV^\gamma_t(\phi)\frac{\xi^2}{2}}.
\end{align*}
We thus concentrate on establishing equation \eqref{eq.crux.central limit theorem}. To this end, note that we can write $A^N_t(\phi) = N^{-1/2} \sum_{i=1}^N U^i_t$,
where the random variables $U^i_t = W^i_t \phi(X^i_t) - \EE_{t-1}\{W^i_t \phi(X^i_t)\} = W^i_t \phi(X^i_t) - \gammahat^N_{t-1}(\forwardQ_t \phi)$ are independent and identically distributed conditionally upon $\F_{t-1}$. Theorem A.3 of \cite{douc2007limit} shows that in order to prove \cref{eq.crux.central limit theorem}, it is enough to prove that for any $\epsilon > 0$ and as $N \plim \infty$, we have
\begin{align} 
\frac{1}{N} \sum_{i=1}^N \Var_{t-1} \BK{ U^i_t } &
\plim
\mu_t(\phi^2) - Z^2_t \pi_{t}(\phi)^2, 
\label{eq.douc.condition}
\\
\frac{1}{N} \sum_{i=1}^N \EE_{t-1}\sqBK{ (U^i_t)^2 \indicator \curBK{ | U^i_t | > N^{1/2} \epsilon } } 
& \plim 
0,
\label{eq.douc.condition.tail}
\end{align} 
where $\indicator (\cdot)$ denotes the indicator function. 
The tail condition \eqref{eq.douc.condition.tail} directly follows from the fact that we consider bounded test functions $\phi \in \B(\XX)$ and that $0 < W^i_t \leq \kappa_g^{t}$ almost surely. We thus focus on proving \cref{eq.douc.condition}. We have $\Var_{t-1} (U^i_t) = \EE_{t-1}\{(W^i_t)^2 \phi^2(X^i_t)\} - \EE_{t-1}\{W^i_t \phi(X^i_t)\}^2$, so \cref{eq.gamma.hat.flow} of the main text, Theorem \ref{thm.consistency}, and the boundedness of $\muhat^N_t(\phi^2)$ together yield that
\begin{align*}
\frac{1}{N} \sum_{i=1}^N \Var_{t-1} \BK{ U^i_t }
& =
\EE_{t-1}\curBK{ \frac{1}{N} \sum_{i=1}^N (W^i_t)^2 \phi^2(X^i_t) }
-
\frac{1}{N} \sum_{i=1}^N \EE_{t-1}\curBK{W^i_t \phi(X^i_t) }^2 \\
& = 
\EE_{t-1}\curBK{ \muhat^N_t(\phi^2)} - \curBK{ \gammahat^N_{t-1}(\forwardQ_t \phi) }^2 
\plim 
\mu_t(\phi^2) - \gamma_{t-1}(\forwardQ_t \phi)^2 \\
& = \mu_t(\phi^2) - \gamma_{t}(\phi)^2 
= \mu_t(\phi^2) - Z^2_t \pi_{t}(\phi)^2,
\end{align*}
as desired. This concludes the proof.
\end{proof}

\bibliographystyle{apalike}
\bibliography{references}

\begin{thebibliography}{}

\bibitem[Ahn et~al., 2014]{ahn2014distributed}
Ahn, S., Shahbaba, B., and Welling, M. (2014).
\newblock Distributed stochastic gradient {MCMC}.
\newblock In {\em International Conference on Machine Learning}, pages
  1044--1052.

\bibitem[Alon, 1986]{alon1986eigenvalues}
Alon, N. (1986).
\newblock Eigenvalues and expanders.
\newblock {\em Combinatorica}, 6(2):83--96.

\bibitem[Andrieu et~al., 2010]{andrieu2010particle}
Andrieu, C., Doucet, A., and Holenstein, R. (2010).
\newblock Particle {M}arkov chain {M}onte {C}arlo methods.
\newblock {\em Journal of the Royal Statistical Society: Series B (Statistical
  Methodology)}, 72(3):269--342.

\bibitem[Andrieu and Roberts, 2009]{andrieu2009pseudo}
Andrieu, C. and Roberts, G.~O. (2009).
\newblock The pseudo-marginal approach for efficient {M}onte {C}arlo
  computations.
\newblock {\em The Annals of Statistics}, pages 697--725.

\bibitem[Beskos et~al., 2014]{beskos2014stability}
Beskos, A., Crisan, D., and Jasra, A. (2014).
\newblock On the stability of sequential {M}onte {C}arlo methods in high
  dimensions.
\newblock {\em The Annals of Applied Probability}, 24(4):1396--1445.

\bibitem[Bolic et~al., 2005]{bolic2005resampling}
Bolic, M., Djuric, P.~M., and Hong, S. (2005).
\newblock Resampling algorithms and architectures for distributed particle
  filters.
\newblock {\em IEEE Transactions on Signal Processing}, 53(7):2442--2450.

\bibitem[Chan and Lai, 2013]{chan2013general}
Chan, H.~P. and Lai, T.~L. (2013).
\newblock A general theory of particle filters in hidden {M}arkov models and
  some applications.
\newblock {\em The Annals of Statistics}, 41(6):2877--2904.

\bibitem[Chopin, 2004]{chopin2004central}
Chopin, N. (2004).
\newblock Central limit theorem for sequential {M}onte {C}arlo methods and its
  application to {B}ayesian inference.
\newblock {\em The Annals of Statistics}, 32(6):2385--2411.

\bibitem[Chopin et~al., 2013]{chopin2013smc2}
Chopin, N., Jacob, P.~E., and Papaspiliopoulos, O. (2013).
\newblock S{MC}2: an efficient algorithm for sequential analysis of state space
  models.
\newblock {\em Journal of the Royal Statistical Society: Series B (Statistical
  Methodology)}, 75(3):397--426.

\bibitem[Del~Moral, 1996]{del1996non}
Del~Moral, P. (1996).
\newblock Non-linear filtering: interacting particle resolution.
\newblock {\em Markov Processes and Related Fields}, 2(4):555--581.

\bibitem[Del~Moral, 2004]{del2004feynman}
Del~Moral, P. (2004).
\newblock {\em Feynman-{K}ac Formulae: Genealogical and interacting particle
  systems with applications, Probability and its applications}.
\newblock Springer.

\bibitem[Del~Moral and Guionnet, 2001]{del2001stability}
Del~Moral, P. and Guionnet, A. (2001).
\newblock On the stability of interacting processes with applications to
  filtering and genetic algorithms.
\newblock {\em Annales de l'Institut Henri Poincare (B) Probability and
  Statistics}, 37(2):155--194.

\bibitem[Del~Moral et~al., 2017]{del2017convergence}
Del~Moral, P., Moulines, E., Olsson, J., and Verg{\'e}, C. (2017).
\newblock Convergence properties of weighted particle islands with application
  to the double bootstrap algorithm.
\newblock {\em Stochastic Systems}, 6(2):367--419.

\bibitem[Douc and Moulines, 2007]{douc2007limit}
Douc, R. and Moulines, E. (2007).
\newblock Limit theorems for weighted samples with applications to sequential
  {M}onte {C}arlo methods.
\newblock In {\em ESAIM: Proceedings}, volume~19, pages 101--107. EDP Sciences.

\bibitem[Durbin and Koopman, 2012]{durbin2012time}
Durbin, J. and Koopman, S.~J. (2012).
\newblock {\em Time Series Analysis by State Space Methods}.
\newblock Number~38 in Oxford Statistical Science Series. Oxford University
  Press.

\bibitem[Fearnhead and Taylor, 2013]{fearnhead2013adaptive}
Fearnhead, P. and Taylor, B.~M. (2013).
\newblock An adaptive sequential {M}onte {C}arlo sampler.
\newblock {\em Bayesian Analysis}, 8(2):411--438.

\bibitem[Friedman, 2008]{friedman2008proof}
Friedman, J. (2008).
\newblock {\em A proof of {A}lon's second eigenvalue conjecture and related
  problems}.
\newblock American Mathematical Society.

\bibitem[Gordon et~al., 1993]{gordon1993novel}
Gordon, N.~J., Salmond, D.~J., and Smith, A.~F. (1993).
\newblock Novel approach to nonlinear/non-{G}aussian {B}ayesian state
  estimation.
\newblock {\em Radar and Signal Processing, IEE Proceedings F},
  140(2):107--113.

\bibitem[Grecian et~al., 2018]{grecian2018understanding}
Grecian, W.~J., Lane, J.~V., Michelot, T., Wade, H.~M., and Hamer, K.~C.
  (2018).
\newblock Understanding the ontogeny of foraging behaviour: insights from
  combining marine predator bio-logging with satellite-derived oceanography in
  hidden {M}arkov models.
\newblock {\em Journal of the Royal Society Interface}, 15(143):20180084.

\bibitem[Hagberg et~al., 2008]{hagberg-2008-exploring}
Hagberg, A.~A., Schult, D.~A., and Swart, P.~J. (2008).
\newblock Exploring network structure, dynamics, and function using {NetworkX}.
\newblock In {\em Proceedings of the 7th Python in Science Conference
  (SciPy2008)}, pages 11--15, Pasadena, CA USA.

\bibitem[Heine and Whiteley, 2017]{heine2017fluctuations}
Heine, K. and Whiteley, N. (2017).
\newblock Fluctuations, stability and instability of a distributed particle
  filter with local exchange.
\newblock {\em Stochastic Processes and their Applications}, 127(8):2508--2541.

\bibitem[Heine et~al., 2020]{heine2020parallelizing}
Heine, K., Whiteley, N., and Cemgil, A.~T. (2020).
\newblock Parallelizing particle filters with butterfly interactions.
\newblock {\em Scandinavian Journal of Statistics}, 47(2):361--396.

\bibitem[Heng and Jacob, 2019]{heng2019unbiased}
Heng, J. and Jacob, P.~E. (2019).
\newblock Unbiased {H}amiltonian {M}onte {C}arlo with couplings.
\newblock {\em Biometrika}, 106(2):287--302.

\bibitem[Ingle et~al., 2015]{ingle2015ultrasonic}
Ingle, A.~N., Ma, C., and Varghese, T. (2015).
\newblock Ultrasonic tracking of shear waves using a particle filter.
\newblock {\em Medical Physics}, 42(11):6711--6724.

\bibitem[Kang et~al., 2018]{kang2018opinion}
Kang, M., Ahn, J., and Lee, K. (2018).
\newblock Opinion mining using ensemble text hidden {M}arkov models for text
  classification.
\newblock {\em Expert Systems with Applications}, 94:218--227.

\bibitem[Lee and Whiteley, 2016]{lee2016forest}
Lee, A. and Whiteley, N. (2016).
\newblock Forest resampling for distributed sequential {M}onte {C}arlo.
\newblock {\em Statistical Analysis and Data Mining: The ASA Data Science
  Journal}, 9(4):230--248.

\bibitem[Lee and Whiteley, 2018]{lee2018variance}
Lee, A. and Whiteley, N. (2018).
\newblock Variance estimation in the particle filter.
\newblock {\em Biometrika}, 105(3):609--625.

\bibitem[Li et~al., 2017]{li2017simple}
Li, C., Srivastava, S., and Dunson, D.~B. (2017).
\newblock Simple, scalable and accurate posterior interval estimation.
\newblock {\em Biometrika}, 104(3):665--680.

\bibitem[Liu and Chen, 1995]{liu1995blind}
Liu, J.~S. and Chen, R. (1995).
\newblock Blind deconvolution via sequential imputations.
\newblock {\em Journal of the American Statistical Association},
  90(430):567--576.

\bibitem[Lorenz, 1963]{lorenz1963deterministic}
Lorenz, E.~N. (1963).
\newblock Deterministic nonperiodic flow.
\newblock {\em Journal of the Atmospheric Sciences}, 20(2):130--141.

\bibitem[Miao et~al., 2011]{miao2011algorithm}
Miao, L., Zhang, J.~J., Chakrabarti, C., and Papandreou-Suppappola, A. (2011).
\newblock Algorithm and parallel implementation of particle filtering and its
  use in waveform-agile sensing.
\newblock {\em Journal of Signal Processing Systems}, 65(2):211--227.

\bibitem[Michelot et~al., 2016]{michelot2016movehmm}
Michelot, T., Langrock, R., and Patterson, T.~A. (2016).
\newblock move{HMM}: an {R} package for the statistical modelling of animal
  movement data using hidden {M}arkov models.
\newblock {\em Methods in Ecology and Evolution}, 7(11):1308--1315.

\bibitem[M{\'\i}guez, 2014]{miguez2014uniform}
M{\'\i}guez, J. (2014).
\newblock On the uniform asymptotic convergence of a distributed particle
  filter.
\newblock In {\em 2014 IEEE 8th Sensor Array and Multichannel Signal Processing
  Workshop (SAM)}, pages 241--244. IEEE.

\bibitem[M{\'\i}guez and V{\'a}zquez, 2016]{miguez2016proof}
M{\'\i}guez, J. and V{\'a}zquez, M.~A. (2016).
\newblock A proof of uniform convergence over time for a distributed particle
  filter.
\newblock {\em Signal Processing}, 122:152--163.

\bibitem[Murray, 2012]{murray2012gpu}
Murray, L. (2012).
\newblock {GPU} acceleration of the particle filter: the {M}etropolis
  resampler.
\newblock {\em arXiv preprint arXiv:1202.6163}.

\bibitem[Murray et~al., 2016]{murray2016parallel}
Murray, L.~M., Lee, A., and Jacob, P.~E. (2016).
\newblock Parallel resampling in the particle filter.
\newblock {\em Journal of Computational and Graphical Statistics},
  25(3):789--805.

\bibitem[Nystrup et~al., 2017]{nystrup2017long}
Nystrup, P., Madsen, H., and Lindstr{\"o}m, E. (2017).
\newblock Long memory of financial time series and hidden {M}arkov models with
  time-varying parameters.
\newblock {\em Journal of Forecasting}, 36(8):989--1002.

\bibitem[Ou et~al., 2021]{ou2021scalable}
Ou, R., Sen, D., and Dunson, D. (2021).
\newblock Scalable {B}ayesian inference for time series via divide-and-conquer.
\newblock {\em arXiv preprint arXiv:2106.11043}.

\bibitem[Qiao et~al., 2017]{qiao2017predicting}
Qiao, F., Li, P., Zhang, X., Ding, Z., Cheng, J., and Wang, H. (2017).
\newblock Predicting social unrest events with hidden {M}arkov models using
  {GDELT}.
\newblock {\em Discrete Dynamics in Nature and Society}, 2017.

\bibitem[Rabiner and Juang, 1986]{rabiner1986introduction}
Rabiner, L. and Juang, B. (1986).
\newblock An introduction to hidden {M}arkov models.
\newblock {\em IEEE ASSP Magazine}, 3(1):4--16.

\bibitem[Scott et~al., 2016]{scott2016bayes}
Scott, S.~L., Blocker, A.~W., Bonassi, F.~V., Chipman, H.~A., George, E.~I.,
  and McCulloch, R.~E. (2016).
\newblock Bayes and big data: The consensus {M}onte {C}arlo algorithm.
\newblock {\em International Journal of Management Science and Engineering
  Management}, 11(2):78--88.

\bibitem[Steger and Wormald, 1999]{steger1999generating}
Steger, A. and Wormald, N.~C. (1999).
\newblock Generating random regular graphs quickly.
\newblock {\em Combinatorics, Probability and Computing}, 8(4):377--396.

\bibitem[Verg{\'e} et~al., 2015]{verge2015parallel}
Verg{\'e}, C., Dubarry, C., Del~Moral, P., and Moulines, E. (2015).
\newblock On parallel implementation of sequential {M}onte {C}arlo methods: the
  island particle model.
\newblock {\em Statistics and Computing}, 25(2):243--260.

\bibitem[Whiteley et~al., 2016]{whiteley2016role}
Whiteley, N., Lee, A., and Heine, K. (2016).
\newblock On the role of interaction in sequential {M}onte {C}arlo algorithms.
\newblock {\em Bernoulli}, 22(1):494--529.

\bibitem[Zhang et~al., 2020]{zhang2020performance}
Zhang, X., Zhao, L., Zhong, W., and Gu, F. (2020).
\newblock Performance analysis of resampling algorithms of parallel/distributed
  particle filters.
\newblock {\em IEEE Access}, 9:4711--4725.

\end{thebibliography}

\end{document}